\documentclass{article}

\usepackage{amsfonts}
\usepackage{amsthm}
\usepackage{amsmath}
\usepackage{xcolor}
\usepackage{comment}
\usepackage{enumerate}

\newtheorem{Thm}{Theorem}[section]
\newtheorem{Lem}[Thm]{Lemma}
\newtheorem{Def}[Thm]{Definition}
\newtheorem{Prop}[Thm]{Proposition}
\newtheorem{Cor}[Thm]{Corollary}

\newtheorem*{Rem}{Remark}

\newcommand{\tb}[1]{\frac{\partial}{\partial \theta _{#1}}}

\newcommand{\sphere}{ \partial B}
\newcommand{\intsphere}{ \int_{\partial B}}

\usepackage{tensor}

\def\sideremark#1{\ifvmode\leavevmode\fi\vadjust{\vbox to0pt{\vss
 \hbox to 0pt{\hskip\hsize\hskip1em
 \vbox{\hsize3cm\tiny\raggedright\pretolerance10000
 \noindent #1\hfill}\hss}\vbox to8pt{\vfil}\vss}}}
\newcommand{\edz}[1]{\sideremark{#1}}

\numberwithin{equation}{section}

\title{Stability of quermassintegral inequalities along inverse curvature flows}
\author{Caroline VanBlargan, Yi Wang }
\date{August 2022}

\begin{document}

\maketitle
\begin{abstract}
In this paper, we consider the stability of quermassintegral inequalities along a inverse curvature flow. We choose a special rescaling of the flow such that the $k$-th quermassintegral is decreasing and the $k-1$-th quermassintegral is preserved. Along this rescaled flow, we prove that the decreasing rate of the $k$-th quermassintegral is faster than the Fraenkel asymmetry of the domain when approaching to the sphere. This leads to the stability inequality of quermassintegral inequalities for nearly spherical sets using the flow method. 
\end{abstract}

\section{Introduction}
\subsection{Quermassintegral inequalities}
For a convex body $\Omega \subset \mathbb{R}^{n+1}$, the \textit{k}-th \textit{quermassintegral} of $\Omega$ is the mixed volume
\begin{align}
    W_k(\Omega) := V(\Omega,...,\Omega,B,...,B),
\end{align}
where $\Omega$ appears in the first $n+1-k$ entries and $B$, which is the unit ball in $\mathbb{R}^{n+1}$, appears in the last $k$ entries. The famous Steiner formula states that the volume of $\Omega + tB$ is a polynomial in $t$. In particular,
\begin{align}
    \text{Vol}(\Omega + tB) = \sum_{k = 0}^{n+1} {n +1 \choose k}W_k(\Omega)t^k.
    \end{align}
Next, denoting $\omega_m$ as the volume of the unit $m$-ball, we set
\begin{align}
    V_k(\Omega) &:= \frac{\omega_k}{\omega_{n+1}}W_{n+1-k}(\Omega).
\end{align}
 Note that $V_{n+1}(A) = \text{Vol}(\Omega)$ and $V_n(\Omega)= \frac{\omega_{n+1}}{(n+1)\omega_n}\text{Area}(\partial \Omega)$. We obtain, as a consequence of the Alexandrov-Fenchel inequalities, the \textit{quermassintegral inequalities}
\begin{align}
   \bigg( \frac{V_{k+1}(\Omega)}{V_{k+1}(B)}\bigg)^{\frac{1}{k+1}} \leq \bigg( \frac{V_{k}(\Omega)}{V_{k}(B)}\bigg)^{\frac{1}{k}}.\label{alek}
\end{align}
When $k=n$, the inequality in \eqref{alek} is simply the classical isoperimetric inequality. For smooth, convex domains and $k\geq 1$, quermassintegrals have the useful integral formula
\begin{align}
\label{mixedvolumeidentity}
    V_{n+1 -k} = \frac{(n+1-k)!(k-1)!}{(n+1)!}\frac{\omega_{n+1-k}}{\omega_{n+1}}\int_{M}  \sigma_{k-1}(L) d\mu,
\end{align}
where  $M:= \partial \Omega$, $L$ is the second fundamental form of $M$, and $\sigma_k(L)$ is the \textit{$k$-th mean curvature} of $M$. The $k$-th mean curvature is the $k$-th elementary symmetric polynomial of the principal curvatures. 
The inequalities in \eqref{alek} equivalently state that for convex domains
\begin{align}
    \bigg (\int_M \sigma_{k-1}(L) d\mu\bigg )^{\frac{1}{n-k+1}}
    \leq
    C(n,k) \bigg (\int_M\sigma_{k}(L) d\mu \bigg )^{\frac{1}{n-k}},
    \label{aleksigma}
\end{align}
where $C(n,k)$ is the constant that gives equality in the case where $M$ is a sphere. More generally, for any $-1\leq j < k$, we have the $(k,j)$\textit{-quermassintegral inequality},
\begin{align}
    \bigg (\int_M \sigma_{j}(L) d\mu\bigg )^{\frac{1}{n-j}}
    \leq
    C(n,k,j) \bigg (\int_M\sigma_{k}(L) d\mu \bigg )^{\frac{1}{n-k}},
    \label{aleksigma1}
\end{align}
where again $C(n,k,j)$ gives equality in the case of the sphere.

 Much of the previous work to establish \eqref{aleksigma1} relies heavily on working with convex domains. There has been work extending \eqref{aleksigma1} to non-convex domains, known as \textit{k-convex domains}, where $\sigma_j(L) \geq 0$ for $1 \leq j \leq k$ with various conditions.

The main goal of this paper is to 
%In Section \ref{sect:flow}, we introduce some results on curvature flows that have been used to study the quermassintegral inequalities in the non-convex case. We further 
study stability of the $(k,j)$\textit{-quermassintegral inequality} \eqref{aleksigma1} along some curvature flows which we will specify later. We consider surfaces $M(t)$, which at each time $t$ is given by an embedding $X: S^n \rightarrow \mathbb{R}^{n+1}$ and satisfies 
\begin{align}
\label{flow_form}
    X_t = G \nu.
\end{align}
$G$ is a symmetric function of the principal curvatures and $\nu$ is the outward pointing vector. Notably, to prove the $(k,k-1)$-quermassintegral inequality for $k$-convex starshaped domains, Guan and Li in \cite{MR2522433} studied the flow
\begin{align}
\label{inverseflowstatement0}
    X_t &= \frac{\sigma_{k-1}(L)}{\sigma_{k}(L)}\nu.
\end{align}
Urbas in \cite{MR1082861} and Gerhardt in \cite{MR1064876} show that the solution exists for all $t \geq 0$ with any initial surface $M(0)$ that is smooth, strictly $k$-convex, and starshaped (and these conditions are preserved for $M(t)$ for all $t\geq 0$). Furthermore, they apply a rescaling to obtain the surfaces $\{\tilde{M}(t)\}$, which they showed converges to a sphere. Additionally, Guan and Li proved that $\frac{d}{dt}\int_{\tilde{M}(t)} \sigma_{k}(\tilde{L}) d \mu_t \leq 0$  and $\frac{d}{dt}\int_{\tilde{M}(t)} \sigma_{k-1}(\tilde{L}) d \mu_t = 0$. These results, combined with an approximation argument to include non-strictly $k$-convex domains, prove the $(k,k-1)$-quermassintegral inequalities. Because of this, to derive stability of the $(k,k-1)$\textit{-quermassintegral inequality} \eqref{aleksigma1}, our most natural flow to consider is the inverse curvature flow \eqref{inverseflowstatement0}.

Later in Section \ref{sect:flowk} we consider a variant flow which preserves volume: 
\begin{align}
    X_t &=  (-\sigma_k(L) + h(t))\nu.
\end{align}
and derive stability of the $(k,-1)$-quermassintegral inequality. The results on this flow are a bit more restrictive. In \cite{MR2647128}, Cabezas-Rivas and Sinestrari found solutions existing for certain convex domains that satisfied a pinching condition. Just like the argument in \cite{MR2522433}, they prove a monotonicity result of the quermassintegral, which provides a proof of the $(k,-1)$ quermassintegral inequalities for this restricted class of surfaces. 

We remark here that in \cite{MR3107515}, Chang and Wang  were able to show \eqref{aleksigma} without the requirement of a starshaped domain, but with the added assumption of having $(k+1)$-convexity instead of just $k$-convexity, and the constant $C(n,k)$ is non-optimal. They proved this using optimal transport methods. See also \cite{MR3291634}, \cite{MR2831436}, and \cite{MR3247383}.

\subsection{Brief discussion of quantitative isoperimetric inequality}

To analyze stability of isoperimetric inequality, one studies the \textit{isoperimetric deficit} $\delta(\Omega)$ of a domain $\Omega \subseteq \mathbb{R}^{n+1}$, defined as
\begin{align}
\label{classical_iso_def}
    \delta(\Omega) := \frac{P(\Omega) - P(B_{\Omega})}{P(B_{\Omega})}.
\end{align}

%Our work to find a quantitative quermassintegral inequality is largely motivated by work done on the stability in the classical isoperimetric inequality. In particular, we are motivated by analysis done on the \textit{isoperimetric deficit} $\delta(\Omega)$ of a domain $\Omega \subseteq \mathbb{R}^{n+1}$, defined as \begin{align} \label{classical_iso_def}    \delta(\Omega) := \frac{P(\Omega) - P(B_{\Omega})}{P(B_{\Omega})}.\end{align} Here $|B_{\Omega}|$ is the volume of $B_{\Omega}$, $B_{\Omega}$ is a ball such that $|\Omega| = |B_{\Omega}|$, and $P( \cdot )$ gives the perimeter of a set. 

The classical isoperimetric inequality is equivalent to $\delta(\Omega) \geq 0$, with equality if and only $\Omega$ is a ball. There has been a lot of work studying quantitative isoperimetric inequalities inspired by the \textit{Bonnesen type inequalities}, which was named by Osserman in \cite{MR519520}. 

This was based off work by Bonnesen, where he studied inequalities in the form \begin{align}  L^2- 4\pi A \geq \lambda(C), \end{align} where $\lambda(C)$ usually denotes some quantity that measures the difference between $\Omega$ to a ball. 
%In this setting, $L$ and $A$ represent the length and area enclosed by a simple closed curve $C$ in $\mathbb{R}^2$. Moreover, $\lambda(C)$ satisfies three conditions: \begin{enumerate}    \item  $\lambda(C) \geq 0$. \item  $\lambda(C) = 0$ precisely when $C$ is a circle.    \item   $\lambda(C)$ measures geometrically how close $C$ is to a circle. \end{enumerate}
Fuglede worked to expand these results to higher dimensions in \cite{MR859955} and \cite{MR942426}, where they proved a quantitative isoperimetric inequality for nearly spherical sets. 
%They used this result to study the stability for convex domains using the spherical deviation. \begin{definition}\label{sphericaldeviation} For a domain $\Omega \subseteq \mathbb{R}^{n+1}$, set    $\tilde{\Omega}:= \frac{\text{Vol}(B)}{\text{Vol}(\Omega)}(\Omega - \text{bar}(\Omega))$,where $\text{bar}(\Omega)$ is the barycenter of $\Omega$ and $B$ is the unit ball. The \textit{spherical deviation} of $\Omega$, $d(\Omega)$, is defined as \begin{align} d(\Omega) := d_H ( \tilde{\Omega}, B ),\end{align} where $d_H(\cdot, \cdot)$ gives the Hausdorff distance between two sets. \end{definition}\noindent Specifically, for a convex domain $\Omega$, they established an inequality in the form $d(\Omega) \leq f(\delta(\Omega))$, for some function $f$.

To establish a quantitative isoperimetric inequality for more general domains, the \textit{Fraenkel asymmetry}, $\alpha(\Omega)$, is a well-studied quantity used as a lower bound. 
\begin{Def}
\label{frankael_def}
Suppose $\Omega \subseteq\mathbb{R}^{n+1}$. The \textit{Fraenkel asymmetry} of $\Omega$ is denoted by $\alpha(\Omega)$, where
\begin{align}
    \alpha(\Omega):= \inf \bigg \{\frac{|\Omega \Delta (x + B_{\Omega})|}{|B_{\Omega}|}, x \in \mathbb{R}^{n+1} \bigg\}.
\end{align}
$B_{\Omega}$ denotes the ball centered at the origin with the same volume as $\Omega$, and $\Delta$ denotes the symmetric difference between two sets.
\end{Def}
Using the Fraenkel asymmetry in the study of stability brings us to the \textit{quantitative isoperimetric inequality}, which asks if there is a fixed $C(n)>0$  such that all Borel sets $\Omega\subseteq \mathbb{R}^{n+1}$ with finite measure satisfy the inequality
\begin{align}
    \delta(\Omega) \geq C(n) \alpha^{m}(\Omega), \label{stability}
\end{align}
for some exponent $m$. The quantitative isoperimetric inequality $ \delta(\Omega) \geq C \alpha^{2}(\Omega)$ was shown for Steiner symmetrical sets in \cite{MR1243100} by  Hall, Hayman, and Weitsman. Later, by using results on the Steiner symmetrical, Hall showed
\eqref{stability} in \cite{MR1166511}, but with a suboptimal exponent of $m=4$.

 Finally, in \cite{MR2456887}, Fusco, Maggi, and Pratelli showed, by using symmetrizations of $\Omega$, that \eqref{stability} holds with optimal exponent $m=2$ for Borel sets $\Omega\subseteq \mathbb{R}^{n+1}$ of finite measure. Figalli, Maggi, and Pratelli in  \cite{MR2672283} proved this optimal result in the more general setting of the anisotropic perimeter.
% \begin{Thm} [\cite{MR2672283}, \cite{MR2456887}]Suppose $n\geq 1$. Then for any Borel set $\Omega\subseteq \mathbb{R}^{n+1}$ of finite measure, \begin{align}     \delta(\Omega) \geq C(n) \alpha^{2}(\Omega), \label{fuscoresult} \end{align} where $C(n) >0$ depends only on $n$. \end{Thm}\noindent This theorem was proved in \cite{MR2456887} by reducing the problem to $n$-symmetric sets and then using the method of Steiner symmetry. To prove the result in \cite{MR2672283}, the authors did not use symmetrization arguments, as done previously. Instead, they applied the Brenier map to employ methods in mass transportation theory. For further reading on the quantitative isoperimetric inequality, see \cite{MR3404715} and \cite{MR2402947}. 

Our study of stability in the quermassintegral inequalities is inspired by work done by Fuglede in \cite{MR942426} and by Cicalese and Leonardi in \cite{MR2980529} on nearly spherical sets.
 \begin{Def}
 \label{nearlysphere}
Suppose $M = \{ (1+u(x))x: x \in \partial B  \}$, where $\partial B$ is the unit sphere in $\mathbb{R}^{n+1}$ and $u:\partial B \rightarrow (-1, \infty)$ is a smooth function on the unit sphere. $M$ is referred to as a \textit{nearly spherical set} when we have suitable, small bounds on $|u|,$ $|\nabla u|$, and $|D^2u|$.
\end{Def}

%Our study of stability in the quermassintegral inequalities is inspired by work done by Fuglede in \cite{MR942426} and by Cicalese and Leonardi in \cite{MR2980529} on nearly spherical sets. \begin{definition} \label{nearlysphere}Suppose $M = \{ (1+u(x))x: x \in \partial B  \}$, where $\partial B$ is the unit sphere in $\mathbb{R}^{n+1}$ and $u:\partial B \rightarrow (-1, \infty)$ is a smooth function on the unit sphere. $M$ is referred to as a \textit{nearly spherical set} when we have suitable, small bounds on $|u|,$ $|\nabla u|$, and $|D^2u|$.\end{definition} \begin{remark} Nearly spherical sets in \cite{MR859955} and \cite{MR2980529} only require bounds on $|u|$ and $|\nabla u|$. However, since we will be working with curvature terms, we will require small bounds on $|D^2u|$ as well. \end{remark}
In \cite{MR2980529}, Cicalese and Leonardi introduced a new method to show \eqref{stability} with optimal exponent $m=2$ for all Borel sets of finite measure. In this paper, they utilized the results in  \cite{MR859955}, which they reformulated by assuming $||u||_{W^{1,\infty}} < \epsilon$ and found
\begin{align}
     \delta(\Omega) \geq  \frac{1 - O(\epsilon)}{2}||u||_{L^2}^2 + \frac{1}{4}|| \nabla u||^2_{L^2}, \label{cicalese_stable}
\end{align}
where $u$ is the function from Definition \ref{nearlysphere}. Note that functions in $O(\epsilon)$ may obtain either positive or negative values. 
It quickly follows that for nearly spherical sets, 
 \begin{align}
 \label{sphericaliso}
     \delta(\Omega) \geq C(1 + O(\epsilon))\alpha^2(\Omega).
 \end{align}
They proved the \textit{Selection Principle}, which provided a new proof of the quantitative isoperimetric inequality by reducing the problem to nearly spherical sets converging to the unit ball. 

It is worth noting that to get \eqref{sphericaliso} from \eqref{cicalese_stable}, only the weaker statement
$ \delta(\Omega) \geq  \frac{1 - O(\epsilon)}{2}||u||_{L^2}^2$ is needed. However, in \cite{MR3216839}, Fusco and Julin showed a stronger result for stability in the isoperimetric problem, where they bound the \textit{asymmetry index} $A(\Omega)$, so that $\delta(\Omega) \geq C A^2(\Omega) $. To do this, they need the full result that $\delta(\Omega) \geq C||u||_{W^{1,2}}^2 $ for nearly spherical sets. Then, as in \cite{MR2980529}, they are able to use the methods of the Selection Principle to reduce the general problem to the results for nearly spherical sets.

In our previous work \cite{MR2441221}, we proved a higher order analogues of the inequality \eqref{cicalese_stable} Cicalese and Leonardi derived in \cite{MR2980529}, and used it to prove the quantitative quermassintegral inequalities.
Let us start by defining higher order isoperimetric deficit function.

%We aim to prove a version of \eqref{cicalese_stable} as it applies to the quermassintegral inequalities. 
First, we define $I_k(\Omega)$ by integrating the $k$-th mean curvature of $M$ for $k \geq 0$,  that is
\begin{align}
     I_k(\Omega) := \int_{M} \sigma_k(L)  \textit{ dA} . 
\end{align}
Also, for $k=-1$ we define
\begin{align}
    I_{- 1}(\Omega) := \text{Vol}(\Omega). 
\end{align}
We are now able to define a natural generalization of the isoperimetric deficit for quermassintegrals.
\begin{Def}
\label{km_deficit}
For $- 1<k \leq n$ and $- 1 \leq m < k$, the \textit{(k,m)-isoperimetric deficit} is denoted by $\delta_{k,m}(\Omega)$, where
\begin{align}
    \delta_{k,m}(\Omega := \frac{I_k(\Omega) - I_k(B_{\Omega,m})}{I_k(B_{\Omega,m})}. 
\end{align}
Here $B_{\Omega,m}$ is the ball centered at the origin where $I_m(B_{\Omega,m})= I_m(\Omega)$. 
\end{Def}

One of our main theorems in \cite{MR2441221} states that
\begin{Thm}(\cite{MR2441221}[Theorem 1.3])
\label{secondthm}
Fix $0 \leq j < k$. Suppose $\Omega= \{ (1+u(\frac{x}{|x|}))x: x \in B  \}\subseteq \mathbb{R}^{n+1}$, where $u \in C^3(\partial B)$,  $I_j(\Omega) = I_j(B)$, and $\text{bar}(\Omega) =0$.
For all $\eta>0$, there exists $\epsilon>0$ such that if  $||u||_{W^{2,\infty}}< \epsilon$, then
\begin{align}
    \delta_{k,j}(\Omega)
    \geq 
\bigg (
\frac{n(n-k)(k-j)}{4(n+1)^2}  - \eta \bigg)\alpha^2(\Omega).
\end{align}
\end{Thm}

 We remark for sufficiently small $||u||_{W^{2,\infty}}$ that $\Omega$ is a convex domain. Then, we already know from the result of Guan and Li, which assumes $\Omega$ is k-convex and starshaped, that $\delta_{k,j}(\Omega) \geq 0$. So, Theorem \ref{secondthm} is establishing a quantitative isoperimetric inequality in this case.

\subsection{Main results}

We now turn our attention specifically to the stability of the $(k,k-1)$-quermassin-tegral inequalities along curvature flows. A key step in \cite{MR2441221} is Proposition 5.1 
%we proved to show our theorems Section 3 
which asserts that for a nearly spherical domain $\Omega$ satisfying $I_{k-1}(\Omega) = I_{k-1}(B)$, \begin{align} \label{toprove}I_k(\Omega) - I_k(B) \geq (1 + O(\epsilon))A(\Omega).   \end{align} 

Our first main result in the current paper studies this inequality along the flow \eqref{inverseflowstatement0}. 
\begin{Thm}
\label{main_inversemeancurv_inequality}
Let $M(t)$ be solution to the flow \eqref{inverseflowstatement0}, and
 $\Tilde{M}(t)$ be the rescaled surface given by 
 \begin{align}
 \label{scaledsurface0}
     \Tilde{X} &= e^{-rt}X,
   \hspace{.15in}
   r = \frac{{n \choose k-1}}{{n \choose k}}.
 \end{align}
Suppose at $t_0$ that $M(t_0)$ is nearly spherical with $||u(t_0)||_{W^{2,\infty}} < \epsilon$  and that the barycenter of $\Tilde{M}(t_0)$ satisfies $|bar(\tilde{M}(t_0))| \leq K \epsilon ||u(t_0)||_{W^{2,2}}^2$ for fixed a $K>0$. Then, for any small $\eta>0$,
\begin{align}
    \frac{d}{dt}(I_k(\tilde{\Omega}(t_0)) - I_k(B))
    \leq 
    (1-\eta)\frac{d}{dt}
    A(\tilde{\Omega}(t_0)),
\end{align}
and the choice of a sufficiently small $\epsilon>0$ depends on $\eta$ and $K$.

Moreover, along any solution to the flow \eqref{inverseflowstatement}  where
$|bar(\tilde{M}(t))| \leq K \epsilon ||u||_{W^{2,2}}^2$ holds for sufficiently large $t$, we have
\begin{align}
   \liminf_{t\rightarrow \infty}  \frac{I_k(\tilde{\Omega}(t)) - I_k(B)}
     {A(\tilde{\Omega}(t))}
   \geq 1.
\end{align}
\end{Thm}
In the following, we will denote \begin{align}\label{A} A(t) &:=A(\Omega(t))=       {n  \choose k}\frac{(n-k)}{2n}\bigg( || u||_{L^2}^2  + \frac{1}{2}||\nabla u||_{L^2}^2  \bigg )\end{align}
for simplicity.

With some additional work, we are able to show 
the following corollary, which provides an alternative proof to the key inequality in \eqref{toprove} in the case where $M$ is $n$-symmetric. 
%\begin{comment}
\begin{Cor}
\label{corollary_inversemean}
Given any $\eta>0$, there is an $\epsilon>0$ such that any smooth $n$-symmetric, nearly spherical set $M$ that, where $M = \partial \Omega$, satisfies the inequality
\begin{align}
    I_k(\Omega) - I_k(B) \geq 
    (1- \eta)A, 
\end{align}
 when $||u||_{C^2}<\epsilon$ and $I_{k-1}(\Omega) = I_{k-1}(B)$.
\end{Cor}
%\end{comment}

%in Corollary \ref{corollary_inversemean} that this provides an alternative proof to the key inequality in \eqref{toprove} in the case where $M$ is $n$-symmetric. Therefore, we view Theorem \ref{main_inversemeancurv_inequality} as a stronger statement to Theorem \ref{secondthm} in this case, as it gives us additional information on the behavior of the derivatives of relevant quantities along the flow. 

We will also study stability of the $(k-1,-1)$-quermassintegral inequality along the volume preserving flow:
%\edz{Either write ``$(k-1,-1)$-quermassintegral inequality", or change the definition of the flow to be $ X_t =  (-\sigma_{k+1}(L) + h(t))\nu.$}
\begin{align}
\label{MeanCurvFlowStatement0}
    X_t &=  (-\sigma_k(L) + h(t))\nu.
\end{align}
The results on this flow are a bit more restrictive. In \cite{MR2647128}, Cabezas-Rivas and Sinestrari found solutions existing for certain convex domains that satisfied a pinching condition. Just like the argument in \cite{MR2522433}, they prove a monotonicity result of the quermassintegral, which provides a proof of the $(k,-1)$ quermassintegral inequalities for this restricted class of surfaces. 

%In \cite{MR3107515}, Chang and Wang  were able to show \eqref{aleksigma} without the requirement of a starshaped domain, but with the added assumption of having $(k+1)$-convexity instead of just $k$-convexity, and the constant $C(n,k)$ is non-optimal. They proved this using optimal transport methods. See also \cite{MR3291634}, \cite{MR2831436}, and \cite{MR3247383}.

Using this flow \eqref{MeanCurvFlowStatement0}, we prove the following theorem.

%We also look at the flow \eqref{MeanCurvFlowStatement0} to prove the following theorem. 

\begin{Thm}
\label{MeanFlowMainThm}
Suppose $M(t)$ is a solution of surfaces to the flow \eqref{MeanCurvFlowStatement0}, and at $t_0$ the surface $M(t_0)$ satisfies, for a fixed $K$, that $|\text{bar}(M(t_0))|^2 \leq K\epsilon||u||_{W^{2,2}}^2$ and $||u(t_0)||_{W^{2,\infty}} < \epsilon$. Then, for any small $\eta>0$
\begin{align}
\label{meanflow_firstineq}
       \frac{d}{dt}\bigg ( I_{k-1}(\Omega(t_0)) - I_{k-1}(B) \bigg )
       \leq
       (1-\eta)\frac{d}{dt}
     \frac{k(n-k+1)}{2n}
    ||u(t_0)||_{L^2}^2,
\end{align}
 where the choice of a sufficiently small $\epsilon>0$ depends on $\eta$ and $K$.
Additionally, if $|\text{bar}(M(t))|^2 \leq K\epsilon||u||_{W^{2,2}}^2$ holds for sufficiently large $t$, then
\begin{align}
\label{liminf_meancurv0}
    \liminf_{t\rightarrow \infty}
    \frac{I_{k-1}(\Omega(t)) - I_{k-1}(B)}{
    ||u||_{L^2}^2} &\geq \frac{k(n-k+1)}{2n}.
\end{align}
\end{Thm}
A few things have prevented us in getting the same results here that we have in Theorem \ref{main_inversemeancurv_inequality}. First, we do not have $\frac{d}{dt}||\nabla u||_{L^2}$ on the right-hand side of \eqref{meanflow_firstineq}. Also, we were not able to use this result to get a quantitative quermassintegral inequality, as we did in a corollary to Theorem \ref{main_inversemeancurv_inequality}. This is because we miss a key lemma stating that if $M(t_0)$ is nearly spherical, then $M(t)$ remains nearly spherical for all $t\geq t_0$. However, this theorem still provides additional information on the stability in the $(k,-1)$-quermassintegral inequalities along the flow.

\section{Preliminaries}
\subsection{The $k$-th mean curvature} \label{sigmaproperty}
For $\lambda = (\lambda_1, ..., \lambda_n) \in \mathbb{R}^n$, we denote $\sigma_k(\lambda)$ as the $k$\textit{-th elementary symmetric polynomial} of $(\lambda_1, ..., \lambda_n)$. That is, for $1 \leq k \leq n$, 

\begin{align}
\sigma_k(\lambda) = \sum_{i_1<i_2<...<i_n} \lambda_{i_1} \lambda_{i_2} \cdot \cdot \cdot \lambda_{i_k},
\end{align}
and 
\begin{align}
\sigma_0(\lambda) = 1.
\end{align}
This leads to a natural generalization of the mean curvature of a surface. 
\begin{Def}
Suppose $\Omega$ is a smooth, bounded domain in $\mathbb{R}^{n+1}$. For $x \in M:=\partial \Omega$, the $k$\textit{-th mean curvature} of $M$ at $x$ is 
$\sigma_{k}(\lambda)$, where $\lambda = (\lambda_1(x), ..., \lambda_n(x))$ are the principal curvatures of $M$ at $x$. 
\end{Def}
Observe that in this definition, $\sigma_1(\lambda)$ is the mean curvature and $\sigma_n(\lambda)$ is the Gaussian curvature. 
When $(\lambda_1,..., \lambda_n)$ are the eigenvalues of a matrix $A = \{A_{j}^i \}$, we denote $\sigma_k(A) = \sigma_k(\lambda)$, which can be equivalently calculated as
\begin{align}
\sigma_k(A) = \frac{1}{k!} \delta_{i_1\cdot\cdot\cdot i_k}^{j_1\cdot\cdot\cdot j_k} A_{j_1}^{i_1}\cdot\cdot\cdot A_{j_k}^{i_k},
\end{align}
 using the Einstein convention to sum over repeated indices.
 
 So, if $L$ is the second fundamental form of $ M$, we can use this expression for $\sigma_k(L)$ to compute the $k$-th mean curvature of $ M$. 
Throughout this paper, we will be working with a family of surfaces where, for $0<j\leq k$, $\sigma_j(L)\geq 0$ at each point. Such surfaces are called \textit{k-convex}.

\begin{Def}
Let $\Omega $ be a domain in $ \mathbb{R}^{n+1}$. Then the hypersurface $M:= \partial \Omega$ is said to be \textit{strictly k-convex} if the principal curvatures $\lambda =(\lambda_1,..., \lambda_n)$ lie in the \textit{G\r{a}rding cone} $\Gamma^{+}_k$, which is defined as
\begin{align}
    \Gamma^{+}_k := \{ \lambda \subseteq 
    \mathbb{R}^n: \sigma_j(\lambda) > 0, 1 \leq j \leq k\}.
\end{align}
\end{Def}
Note that n-convexity is the same as normal convexity. A useful operator related to $\sigma_k$ is the  \textit{Newton transformation tensor} $[T_k]^j_i$.
\begin{Def}
The Newton transformation tensor, $[T_k]^j_i$, of $n\times n$ matrices $\{ A_{1}, ..., A_k\}$ is defined as
\begin{align}
    [T_k]^j_i(A_1,...,A_k):= \frac{1}{k!} \delta_{ii_1\cdot\cdot\cdot i_k}^{jj_1\cdot\cdot\cdot j_k} (A_1)_{j_1}^{i_1}\cdot\cdot\cdot (A_k)_{j_k}^{i_k}.
\end{align}
When $A_1 = A_2 =...=A_k = A$, we denote $[T_k]^j_i(A) = [T_k]^j_i(A,... ,A)$.
\end{Def}
A related operator is  $\Sigma_k$, which  the polarization of $\sigma_k$.
\begin{Def}
Suppose $\{ A_{1}, ..., A_k\}$ is a collection of $n\times n$ matrices. We denote
    \begin{align}
    \Sigma_k(A_1,...,A_k)
    &:= 
    (A_1)_{j}^i[T_{k-1}]^j_i(A_2,...,A_k)
   \nonumber
   \\ 
   &= \frac{1}{(k-1)!} \delta_{i_1\cdot\cdot\cdot i_k}^{j_1\cdot\cdot\cdot j_k} (A_1)^{i_1}_{j_1}\cdot\cdot\cdot (A_k)^{i_k}_{j_k}.
\end{align}
\end{Def}
Two useful identities are
\begin{align}
    \sigma_k(A) = \frac{1}{k} \Sigma_k(A,...,A) = \frac{1}{k} A_{j}^i[T_{k-1}]^j_{i}(A) ,
\end{align}
and
\begin{align}
    \frac{\partial \sigma_k(A)}{\partial A_{j}^i} = \frac{1}{k}[T_{k-1}]^j_i(A).
\end{align}
We will also use the identity 
\begin{align}
    A_{s}^j[T_m]^i_j(A) = \delta^{i}_s\sigma_{m+1}(A) - [T_{m+1}]_s^i(A).
\end{align}

\subsection{Nearly spherical sets}
\label{sphericalprelim}

   The focus of this paper is to establish the $(k,m)$-isoperimetric inequality for \textit{nearly spherical sets}. Our approach is inspired by Cicalese and Leonardi's work in the classical quantitative isoperimetric inequality for nearly spherical sets in \cite{MR2980529}.
That is, we consider a smooth, bounded domain $\Omega$ that is starshaped with respect to the origin, which is enclosed by $M := \partial \Omega$. We write $M = \{ (1+u(x))x: x \in \partial B  \}$, where $u:\partial B \rightarrow \mathbb{R}$ is a smooth function. The set $M$ is referred to as a \textit{nearly spherical set} when there is a suitable, small bound on $||u||_{W^{2,\infty}}$. In this section, we establish some useful formulas for nearly spherical sets.

We write $\mathbb{R}^{n+1}$ in spherical coordinates with the tangent basis $\{\frac{\partial}{\partial \theta_1},\frac{\partial}{\partial \theta_2}, ...,$
$\frac{\partial}{\partial \theta_n}, \frac{\partial}{\partial r} \}$. Denoting $s_{ij}$ as the metric on the sphere, we have   $<\frac{\partial}{\partial \theta_i}, \frac{\partial}{\partial r} > = 0$, $<\frac{\partial}{\partial r}, \frac{\partial}{\partial r} > = 1$, and $<\frac{\partial}{\partial \theta_i}, \frac{\partial}{\partial \theta_j} > = r^2s_{ij}$.
Set $u_i = \frac{\partial u}{\partial \theta_i}$. Then, $\{e_i\}$ forms a tangent basis of $M$ where
\begin{align}
e_i = \frac{\partial}{\partial \theta_i} + u_i\frac{\partial}{\partial r}. 
\end{align}
We find,
\begin{align}
\label{outernormal_formula}
N = \frac{-\sum_{i=1}^n s^{ij}u_i \frac{\partial}{\partial \theta_j}+ (1+u)^2 \frac{\partial}{\partial r} }{(1+u)\sqrt{|\nabla u|^2 + (1+u)^2}} ,
\end{align}
where  $N$ is the outward  unit normal on $M$, and the norm $|\nabla u|$ is taken with respect to the standard metric on $\partial B$.
We compute the metric $g_{ij}$ on $M$ as
\begin{align}
g_{ij} = <e_i,e_j> = (1+u)^2s_{ij} + u_iu_j , 
\end{align}
where $<\cdot, \cdot>$ is the standard Euclidean inner product on $\mathbb{R}^{n+1}$. Setting $g^{ij}$ to be the inverse of $g_{ij}$, we have

\begin{align}
g^{ij}  = \frac{s^{ij}}{(1+u)^2} - \frac{1}{(1+u)^2}\frac{u_ku_ls^{ki}s^{lj}}{|\nabla u|^2 + (1+u)^2}  .
\end{align}
We denote $h_{ij}$ as the second fundamental form on $M$. That is, 
$h_{ij}=-<N,\nabla_{e_i}e_j>$, and we form the shape operator $h^i_j$ by
\begin{align}
h^i_j = g^{ik}h_{kj}.
\end{align}
We now explicitly calculate $h^i_j$. First, note
\begin{align}
   \bigg (\nabla_{\frac{\partial}{\partial \theta_i}}\frac{\partial}{\partial \theta_j}
   \bigg)^k
    &=
    \frac{1}{2}\frac{1}{r^2}s^{kl}
    \bigg(
    \partial_i (r^2 s_{il}) + \partial_j(r^2s_{il}) - \partial_l(r^2s_{ij})
    \bigg)
    \\
    &=
    \frac{1}{2}s^{kl}
    \bigg(
    \partial_i s_{il} + \partial_js_{il} - \partial_ls_{ij}
    \bigg)
    \\
    &= \Gamma^k_{ij},
\end{align}
where $\Gamma^k_{ij}$ refers to the Christoffel symbol on $\mathbb{S}^n$, and \begin{align}
   \bigg (\nabla_{\frac{\partial}{\partial \theta_i}}\frac{\partial}{\partial \theta_j}
   \bigg)^r
    &=-rs_{ij}.
\end{align}
We thus obtain
\begin{itemize}
    \item $\nabla_{\frac{\partial}{\partial \theta_i}}\frac{\partial}{\partial \theta_j} = \Gamma^k_{ij}\frac{\partial}{\partial \theta_k}
    -rs_{ij}\frac{\partial}{\partial r}$.\end{itemize}
Similarly,    
    \begin{itemize}
    \item$ \nabla_{\frac{\partial}{\partial \theta_i}}\frac{\partial}{\partial r} 
    =
    \nabla_{\frac{\partial}{\partial r}}\frac{\partial}{\partial \theta_i}
    =
    \frac{1}{r}\frac{\partial}{\partial \theta_i}$,
       \item$ \nabla_{\frac{\partial}{\partial r}}\frac{\partial}{\partial r} = 0$.
\end{itemize}
Then,
\begin{align}
    \nabla_{e_i} e_j 
      =& \nabla_{\frac{\partial}{\partial \theta_i} + u_i\frac{\partial}{\partial r} }
      \bigg(
      \frac{\partial}{\partial \theta_j} + u_j\frac{\partial}{\partial r}
      \bigg)
      \nonumber 
      \\
      =&  \nabla_{\frac{\partial}{\partial \theta_i} }\frac{\partial}{\partial \theta_j} 
      +  \nabla_{\frac{\partial}{\partial \theta_i} }
      \bigg( u_j\frac{\partial}{\partial r} \bigg)
      + u_i\nabla_{ \frac{\partial}{\partial r} }\frac{\partial}{\partial \theta_j}
       + u_i\nabla_{ \frac{\partial}{\partial r} } 
       \bigg(
       u_j\frac{\partial}{\partial r}
       \bigg)
       \nonumber 
       \\
       =&  \nabla_{\frac{\partial}{\partial \theta_i} }\frac{\partial}{\partial \theta_j} 
      + u_j \nabla_{\frac{\partial}{\partial \theta_i} } \frac{\partial}{\partial r}
      +
      \bigg(\frac{\partial^2}{\partial \theta_i\partial \theta_j}
      u
      \bigg)\frac{\partial}{\partial r} 
      + u_i\nabla_{ \frac{\partial}{\partial r} }\frac{\partial}{\partial \theta_j}
       + u_iu_j\nabla_{ \frac{\partial}{\partial r} }\frac{\partial}{\partial r} \nonumber \\
      & + u_i
       \bigg(\frac{\partial}{\partial r}u_{i}
       \bigg) \frac{\partial}{\partial r} \nonumber 
        \\
       =&  \Gamma^k_{ij}\frac{\partial}{\partial \theta_k}
    -(1+u)s_{ij}\frac{\partial}{\partial r}
      +\frac{1}{(1+u)}(u_j\tb{i} + u_i\tb{j})
      +
      \bigg(\frac{\partial^2}{\partial \theta_i\partial \theta_j}
      u
      \bigg)\frac{\partial}{\partial r}.
\end{align}
So
\begin{align}
 h_{ij}=& -
 \bigg<
 \frac{-s^{pq}u_{p} \frac{\partial}{\partial \theta_q}+ (1+u)^2 \frac{\partial}{\partial r} }{(1+u)\sqrt{|\nabla u|^2 + (1+u)^2}} 
 ,
 \Gamma^k_{ij}\frac{\partial}{\partial \theta_k}
    -(1+u)s_{ij}\frac{\partial}{\partial r}
  \nonumber \\
   &+\frac{1}{(1+u)} \bigg(u_j\tb{i} + u_i\tb{j}
      \bigg) +
      \bigg(\frac{\partial^2}{\partial \theta_i\partial \theta_j}
      u
      \bigg)\frac{\partial}{\partial r}
    \bigg >.
  \end{align}
Thus,
  \begin{align}
h_{ij} =&\frac{1}{\sqrt{|\nabla u|^2 +(1+u)^2}}
 \bigg(
 (1+u)u_k\Gamma^k_{ij}
 + (1+u)^2s_{ij} 
 +
 2u_iu_j   \nonumber \\
     & -  (1+u) 
      \bigg(
 \frac{\partial^2}{\partial \theta_i\partial \theta_j} u
      \bigg)
     \bigg )  \nonumber
       \\
 =&\frac{1}{\sqrt{|\nabla u|^2 +(1+u)^2}} (2u_iu_j
     + (1+u)^2s_{ij} 
      -  (1+u) u_{ij}
    \bigg  ),
\end{align}
where $u_{ij}$ denotes the Hessian of $u$ on $\mathbb{S}^n$.
Set
\begin{align}
D := \sqrt{|\nabla u|^2 +(1+u)^2}.
\end{align}
Then,
\begin{align}
    h^i_{j} 
    =& g^{ik}h_{kj}
     \nonumber
    \\
   =& 
   \bigg(\frac{s^{ik}}{(1+u)^2} - \frac{1}{(1+u)^2}\frac{u_mu_ls^{mi}s^{lk}}{D^2}\bigg)
   \frac{1}{D}
   \bigg(
   2u_ku_j
     + (1+u)^2s_{kj} 
      \nonumber\\
      &-  (1+u) u_{kj}
   \bigg)
    \nonumber
    \\
    =&
    \frac{2u^iu_j}{(1+u)^2D} + \frac{\delta^i_j}{D} - \frac{u^i_{j}}{(1+u)D}
    -\frac{2u^iu_j|\nabla u|^2}{(1+u)^2D^3} -\frac{u^iu_j}{D^3}
    +\frac{u^iu_lu^l_{j}}{(1+u)D^3}.
\end{align}
We observe that
\begin{align}
     \frac{2u^iu_j}{(1+u)^2D}   -\frac{2u^iu_j|\nabla u|^2}{(1+u)^2D^3}
     - \frac{u^iu_j}{D^3}
&= \frac{u^iu_j}{D^3},
\end{align}
which yields
\begin{align}
    h^i_{j} 
    &=
   \frac{\delta^i_j}{D} - \frac{u^i_{j}}{(1+u)D}
    +\frac{u^iu_lu^l_{j}}{(1+u)D^3}
    +\frac{u^iu_j}{D^3}.
\end{align}
Next, note
\begin{align}
\label{area_element}
&\sqrt{ \text{det } g_{ij}  }= (1+u)^n \sqrt{\frac{|\nabla u|^2}{(1+u)^2} + 1}.
\end{align}
Therefore,
\begin{align}
    & \text{Area}(M) = \intsphere (1+u)^n \sqrt{\frac{|\nabla u|^2}{(1+u)^2} + 1} \textit{ dA}.
\end{align}We list a few more relevant formulae below:
\begin{align}
   & |\Omega| = \frac{1}{n+1}\intsphere (1 +u)^{n+1} \textit{dA},
   \label{volume_formula}
    \\
    &|\Omega \Delta B| =    \sum_{k=1}^{n+1} \int_{\sphere} \frac{1}{n+1}{n +1 \choose k}|u|^k \textit{dA},
    \label{symmetric_diff_formula}
    \\
    &\text{bar}(\Omega) = \frac{1}{\text{Area}(\partial B)}\intsphere (1+u)^{n+2} x \textit{dA}
    \label{bar_formula}.
\end{align}

Finally, we consider how to compute $\nabla_j [T_m]^j_i(D^2u)$ for nearly spherical sets. This is particularly useful when applying integration by parts to $I_k(\Omega)$ in Section \ref{sect:flow}. See \cite{MR3107515} for a similar computation.  We compute,
\begin{align}
    \nabla_j [T_m]^{j}_i(D^2u) 
     =&
     \frac{1}{m!}  \nabla_j\delta^{jj_1j_2...j_{m}}_{ii_1i_2...i_{m}} u^{i_1}_{j_1}u^{i_2}_{j_2}\cdot \cdot\cdot u^{i_m}_{j_m}
   \nonumber  \\
     =&\frac{m}{m!} \delta^{jj_1j_2...j_{m}}_{ii_1i_2...i_{m}} 
      (\nabla_ju^{i_1}_{j_1})u^{i_2}_{j_2}\cdot \cdot\cdot u^{i_m}_{j_m}.
    \end{align}
    Note that
  \begin{align}
     \delta^{jj_1j_2...j_{m}}_{ii_1i_2...i_{m}} (\nabla_j u^{i_1}_{j_1})u^{i_2}_{j_2}\cdot \cdot\cdot u^{i_m}_{j_m}
     =&
      -\delta^{j_1jj_2...j_{m}}_{ii_1i_2...i_{m}} (\nabla_ju^{i_1}_{j_1})u^{i_2}_{j_2}\cdot \cdot\cdot u^{i_m}_{j_m}.
  \end{align}
  We obtain
    \begin{align}
       \nabla_j [T_m]^j_i(D^2u)    =&  
       \frac{1}{2(m-1)!}
          \delta^{jj_1j_2...j_{m}}_{ii_1i_2...i_{m}} 
          (\nabla_ju^{i_1}_{j_1} -\nabla_{j_1}u^{i_1}_{j})
          u^{i_2}_{j_2}
          \cdot \cdot\cdot u^{i_m}_{j_m}
     \nonumber  \\
    =&   \frac{1}{2(m-1)!}
    \delta^{jsj_1...j_{m-1}}_{ili_1...i_{m-1}}
    (u_pR^{pl}_{\;\;\;sj})
     u^{i_1}_{j_1}
          \cdot \cdot\cdot u^{i_{m-1}}_{j_{m-1}},
    \label{derivativenewton}
\end{align}
where $R^{pl}_{\;\;\;sj}$ is the curvature tensor on $\Omega$.
On $\mathbb{S}^n$, we know by the Gauss equation,
\begin{align}
    R^{pl}_{\;\;\;sj} &= h^p_{s}h^l_{j} - h^p_{j}h^l_{s}
=\delta^p_{s}\delta^l_{j} - \delta^p_{j}\delta^l_{s}.
\end{align}
Therefore,
 \begin{align}
       \nabla_j [T_m]^j_i(D^2u)   
    =&   
    \frac{1}{(m-1)!}\frac{1}{2}
    u_p(\delta^p_{s}\delta^l_{j} - \delta^p_{j}\delta^l_{s})
    \delta^{jsj_1...j_{m-1}}_{ili_1...i_{m-1}}
     u^{i_1}_{j_1}
          \cdot \cdot\cdot u^{i_{m-1}}_{j_{m-1}}
    \nonumber   
    \\
     =&   
    \frac{1}{(m-1)!}\frac{1}{2}
    \bigg (
    u_s
    \delta^{jsj_1...j_{m-1}}_{iji_1...i_{m-1}}
     u^{i_1}_{j_1}
          \cdot \cdot\cdot u^{i_{m-1}}_{j_{m-1}}
          -
       u_j
    \delta^{jlj_1...j_{m-1}}_{ili_1...i_{m-1}}
     u^{i_1}_{j_1} \nonumber\\
     &\cdots  u^{i_{m-1}}_{j_{m-1}}
          \bigg)
    \nonumber  
    \\
        =&   
    \frac{-1}{(m-1)!}
       u_j
    \delta^{jlj_1...j_{m-1}}_{ili_1...i_{m-1}}
     u^{i_1}_{j_1}
          \cdot \cdot\cdot u^{i_{m-1}}_{j_{m-1}}
\nonumber
          \\
          =&
    -(n-m)
    u_j [T_{m-1}]^j_i(D^2u).
\end{align}

\section{Inverse curvature flow} 
\label{sect:flow}

In \cite{MR2522433}, Guan and Li gave a proof of the quermassintegral inequalities for $k$-convex starshaped domains. They used a special case of the flows studied by Urbas in \cite{MR1082861} and Gerhardt in \cite{MR1064876}, which looked at certain flows involving a function $f$ that is symmetric in its inputs $\lambda= (\lambda_1,...,\lambda_n)$, which are the principal curvatures of a hypersurface in $\mathbb{R}^{n+1}$. We take $\Gamma \subseteq \mathbb{R}^n$ to be any set that is an open, convex cone with its vertex at the origin, which also contains $\Gamma_n =\{ (\lambda_1,...,\lambda_n) \in \mathbb{R}^n:  \lambda_i >0 \text{ for each } 1 \leq i \leq n \}$. The function $f$ is assumed to satisfy the following properties:
\begin{align}
\label{first_property}
    f\in C^{\infty}(\Gamma) \cap C^0(\overline{\Gamma}) ;
    \\
    \frac{\partial f}{\partial \lambda_i} >0 \text{ for } \lambda \in \Gamma;
    \\
    \frac{\partial^2 f}{\partial \lambda_i\partial\lambda_j} \leq0;
    \\
    \label{last_property}
    f \equiv 0 \text{ on }  \partial \Gamma.
\end{align}

We now examine the flow of surfaces $M(t)$, which at each time $t$ is an embedding given by $X: S^n \rightarrow \mathbb{R}^{n+1}$, and satisfies the equation 
\begin{align}
\label{statement_inverse_flow}
    X_t &= \frac{1}{f}\nu,
\end{align}where $\nu$ is the outer unit normal vector on $M(t)$. The following theorem, proved by Urbas and Gerhardt, discusses the existence of starshaped solutions.
\begin{Thm}(\cite{MR1082861} and \cite{MR1064876}]
Suppose $M(0)$ is a smooth surface that is starshaped with respect to a point $p$ and with principal curvatures in $\Gamma$. There is a unique smooth solution of surfaces $\{M(t)\}$  to the flow \eqref{statement_inverse_flow} for $t\in [0, \infty)$, and $M(t)$ remains starshaped with respect to $p$.

Furthermore, when rescaling the surfaces to $\{\tilde{M}(t)\}$ so that $\tilde{X} = e^{-\beta t}X$, where $\beta = f(I)$, $\{\tilde{M}(t)\}$ converges to a sphere exponentially fast.
\end{Thm}
\begin{Rem}
 From here on, we assume the surfaces $\{M(t)\}$ are starshaped with respect to the origin. We denote $1+u(t)$ as the radius of $\Tilde{M}(t)$. As a result of the exponential convergence of $\Tilde{M}(t)$ to the sphere, for each $k \geq 0$ there exist $K,\gamma >0$ where $||D^ku(t)||_{\infty} < K e^{-\gamma t}$.
\end{Rem}

To study the quermassintegral inequalities, Guan and Li in \cite{MR2522433} specifically look at the case where $f = \frac{\sigma_k(L)}{\sigma_{k-1}(L)}$. That is, 
\begin{align}
    X_t &= \frac{\sigma_{k-1}(L)}{\sigma_{k}(L)}\nu.
\end{align}
We state some useful derivatives along the flow \eqref{statement_inverse_flow} from \cite{MR2522433} in the following proposition.
\begin{Prop}
\label{flow_derivatives_prop}
(\cite{MR2522433}[Proposition 4])
%(Guan and Li, \cite{MR2522433} Proposition 4)
\begin{enumerate}
    \item $\partial_t g_{ij} = \frac{2}{f}L_{ij}$
    \item $\partial_t \nu = -\nabla \frac{1}{f}$
     \item $\partial_t d\mu_g = \frac{1}{f}\sigma_1(L)d\mu_g $
      \item $\partial_t L_{ij} = -\nabla_i\nabla_j\frac{1}{f} +\frac{(L^2)_{ij}}{f}$
       \item $\partial_t L^i_j = -\nabla^i\nabla_j\frac{1}{f} -\frac{(L^2)^i_{j}}{f}$
         \item $\partial_t \sigma_m(L) = -\nabla_j([T_{m-1}]^i_j(L) \nabla_i \frac{1}{f})
         -\frac{1}{f}
         [T_{m-1}]^{i}_j(L) (L^2)^j_i$,
         %\Sigma_{m-1,1}(h^i_j,(h^2)^i_{j})$,
\end{enumerate}
where $(L^2)^j_{i}= g^{kl}g^{jp}L_{pk}L_{li}$.

\end{Prop}
%\edz{1. $\frac{1}{m}$ factors need to be removed; \\2. $\Sigma_{m-1,1}(h^i_j,(h^2)^i_j)$ is not defined. Use $T_{m-1}$ notation instead. Be consistent with notation on Page 20 after making the notation change.\\3. All $h$'s should be $L$'s.}

In \cite{MR2522433}, the authors show 
\begin{align}
  \label{differentiate_sigma_along_flow}
 \frac{d}{dt} \int_{M(t)} \sigma_k(L) d\mu_t
    &= (k+1)  \int_{M(t)} \frac{\sigma_{k+1}(L)\sigma_{k-1}(L)}{\sigma_k(L)} d\mu_t.
\end{align}
Their argument extends to general curvature flow in the form $X_t = \frac{1}{f}\nu$ and derive:
\begin{Prop}
(\cite{MR2522433}[Lemma 5])
Along the flow \eqref{statement_inverse_flow}, it holds that
\begin{align}
\frac{d}{dt}\int_{M(t)} \sigma_m(L) d\mu_t
    &=
          (m+1)  \int_{M(t)} \frac{\sigma_{m+1}(L)}{f} d\mu_t.
\end{align}
\end{Prop}

For our computations, we find it convenient to consider a normalization $\tilde{\Omega} (t)$ that fixes $I_{k-1}(\tilde{\Omega} (t)) = I_{k-1}(\tilde{\Omega} (0))$. 
Note that this normalization is different from what is adopted by \cite{MR2522433} in which $I_{k}(\tilde{\Omega} (t)) = I_{k}(\tilde{\Omega} (0))$, but as $I_{k-1}(\tilde{\Omega}(t))$ remains constant and $I_{k}(\tilde{\Omega}(t))$ decreases along the flow, we make the same conclusion as in \cite{MR2522433} to prove the $(k,k-1)$-quermassintegral inequalities. Our choice of normalization is by setting $r(t) \equiv \frac{{n \choose k-1}}{{n \choose k}}$ and $\tilde{X} = e^{-\int_0^tr(s)ds}X$. By
 similar computations, we apply formula \eqref{differentiate_sigma_along_flow} to compute
\begin{align}
\label{diff_scaled_along_inverse_flow}
    \frac{d}{dt} \int_{\Tilde{M}} \sigma_m(\Tilde{L}) d\Tilde{\mu_t}
    =&
    \frac{d}{dt} 
    \bigg( e^{-(n-m) \int_0^t r(s) ds} 
    \int_{{M}} \sigma_m({L}) d\mu_t
    \bigg )
    \nonumber
    \\
    =&
   (m+1)e^{-(n-m) \int_0^t r(s) ds}
   \bigg( 
   \int_M \frac{\sigma_{m+1}(L)}{f} d\mu_t \nonumber\\
   &-\frac{(n-m)}{m+1}r(t)\int_{{M}} \sigma_m({L}) d\mu_t
   \bigg).
  %=&(m+1) e^{-(n-m) \int_0^t r(s) ds}  \bigg( \int_{\Tilde{M}}   \frac{\sigma_{m+1}(\Tilde{L})}{\Tilde{f}} d\mu_t-\frac{{n \choose m +1}}{{n \choose m}}r(t)\int_{{\Tilde{M}}} \sigma_m({\Tilde{L}})d\mu_t\bigg).
\end{align}

%\edz{1. $\tilde{\mu_t}$ instead of $\mu_t$ when proper; 2. we cannot change back to $\tilde{L}$ from the second line to the 3rd line of (3.9), but the computation works if we keep the exponential $r$ term, and do not change back to $\tilde{L}$. We need to rewrite (3.9). 3. The same problem in (3.10), we should use $L$, not $\tilde{L}$.4. The same problem in (3.11), and (3.12).}

Next, using our choice of $r(t) = \frac{{n \choose k-1}}{{n \choose k}}$ for the rescaling and substituting in $\frac{1}{f} = \frac{\sigma_{k-1}(L)}{\sigma_{k}(L)}$, we have
\begin{align}
    \frac{d}{dt} \int_{\Tilde{M}} \sigma_m(\Tilde{L}) d\Tilde{\mu_t}
       =&
       (m+1)
   \bigg( 
   \int_{M}
   \frac{\sigma_{m+1}(L)\sigma_{k-1}(L)}{\sigma_{k}(L)} d\mu_t \nonumber\\
  & -\frac{{n \choose m +1}}{{n \choose m}}\frac{{n \choose k-1}}{{n \choose k}}
   \int_{{M}} \sigma_m({L})
   d\mu_t
   \bigg ).
\end{align}
Applying the Newton-Maclaurin inequalities (as in \cite{MR2522433}), which states $$\sigma_{k+1}(L)\sigma_{k-1}(L) \leq \frac{{n \choose k +1}{n \choose k-1}}{{n \choose k}^2}\sigma_{k}^2(L), $$ we find
\begin{align}
\label{deriv_k}
    &\frac{d}{dt} \int_{\Tilde{M}} \sigma_k(\Tilde{L}) d\Tilde{\mu_t} \nonumber\\
       =&
       (k+1)
   \bigg( 
   \int_{M}
   \frac{\sigma_{k+1}(L)\sigma_{k-1}(L)}{\sigma_{k}(L)} d\mu_t
   -\frac{{n \choose k +1}{n \choose k-1}}{{n \choose k}^2}
   \int_{{M}} \sigma_k({L})
   d\mu_t
   \bigg )
   \leq 0,
\end{align}
and
\begin{align}
\label{deriv_kminus1}
    \frac{d}{dt} \int_{\Tilde{M}} \sigma_{k-1}(\Tilde{L}) d\Tilde{\mu_t}
     &=
     k
   \bigg( 
   \int_{M} 
   \frac{\sigma_{k}(L)\sigma_{k-1}(L)}{\sigma_{k}(L)} d\mu_t
   -
   \int_{M} \sigma_{k-1}({L})
   d\mu_t
   \bigg )
   =0.
   \end{align}

In the next section, we will study the curvature flows to prove quantitative quermassintegral inequalities. In doing so, we need to study the quantities $\frac{d}{dt}||u||_{L^2}^2$ and $\frac{d}{dt}||\nabla u||_{L^2}^2$. For this, we need to understand how to convert the equation $X_t = \frac{1}{f}\nu$ into an equation involving the radius of the surface, instead of $X$. Both Urbas and Gerhardt derive a formula for differentiation of the radial function $w$, which we state in the following lemma. We state the formula in a slightly more general form $X_t = G \nu$. %(and $G$ does not need to satisfy the properties of $f$ stated at the beginning of this section).%This way, we can apply it to other flows, as well. 

%\edz{Should it be "$1/G= f$ does need to satisfy the properties of $f$ stated at the beginning of this section"?}

\begin{Lem}[See \cite{MR1082861} and \cite{MR772132}]
\label{deriv_of_radius}
Suppose $\{M(t)\}$ is a collection of smooth hypersurfaces in $\mathbb{R}^{n+1}$ such that each $M(t)$ is starshaped with respect to the origin and satisfies
\begin{align}
    X_t &= G\nu. 
\end{align}
If we write $M(t)$ in spherical coordinates with coefficients $(\theta_1,...,\theta_n, w)$, where $w$ is the radius, we have
\begin{align}
    w_t &= G\sqrt{1 + \frac{|\nabla w|^2}{w^2}}.
\end{align}
\end{Lem}
%\edz{Do we need to give a proof of Lemma 3.6 here or should we quote \cite{MR1082861} and \cite{MR772132} directly?}

\begin{proof}
%Note, we are using $w_t$ and $\frac{\partial w}{\partial t}$ interchangeably. 
The outer unit normal at a point on $M(t)$ in spherical coordinates is given by
\begin{align}
\nu = \frac{-\sum_{i=1}^n s^{ij}\frac{w_i}{w^2} \frac{\partial}{\partial \theta_j}+  \frac{\partial}{\partial r} }{P} ,
\end{align}
where $P = \sqrt{1 + \frac{|\nabla w|^2}{w^2}}$. For a solution to the flow $X_t = G\nu$, we can write
$    X(\theta, t) = w(\phi(\theta,t),t)\phi(\theta,t)$,
where $w$ is the radial function of $M(t)$ at $\phi(\theta,t)$, and $\phi (\cdot, t): \mathbb{S}^n \rightarrow \mathbb{S}^n$ is a suitable diffeomorphism of $S^n$. Hence
\begin{align}
    X_t
    &= 
    (w_m \frac{\partial \phi^m}{dt} + \frac{\partial w}{\partial t})\phi(\theta) 
    + w \frac{\partial \phi}{\partial t},
\end{align}
and
\begin{align}
    G &=  <X_t, \nu>
    \nonumber
    \\
    &= \frac{1}{P} <   (w_m \frac{\partial \phi^m}{\partial t} + \frac{\partial w}{\partial t})\phi(\theta) + w \frac{\partial \phi}{\partial t}, -\sum_{i=1}^n s^{ij}\frac{w_i}{w^2} \frac{\partial}{\partial \theta_j}+  \frac{\partial}{\partial r}>
     \nonumber
    \\
    &=  \frac{1}{P} <   (w_m \frac{\partial \phi^m}{\partial t} + \frac{\partial w}{\partial t})\phi(\theta),  \frac{\partial}{\partial r}>
    +
     \frac{1}{P} <   w \frac{\partial \phi}{\partial t}, -\sum_{i=1}^n s^{ij}\frac{w_i}{w^2} \frac{\partial}{\partial \theta_j}>
       \nonumber
    \\
  &=  \frac{1}{P} <   (w_m \frac{\partial \phi^m}{\partial t} + \frac{\partial w}{\partial t})\frac{\partial}{\partial r},  \frac{\partial}{\partial r}>
    +
     \frac{1}{P} <   w \frac{\partial \phi^m}{\partial t}\frac{\partial}{\partial \theta_m}, -\sum_{i=1}^n s^{ij}\frac{w_i}{w^2} \frac{\partial}{\partial \theta_j}>
       \nonumber
    \\
    &= \frac{1}{P}
    \bigg (
    w_m \frac{\partial \phi^m}{\partial t} + \frac{\partial w}{\partial t}
    -
    w\frac{\partial \phi^m}{\partial t}s^{ij}\frac{w_i}{w^2}s_{jm}w^2
    \bigg )
           \nonumber
    \\
    &= \frac{1}{P}
    \bigg (
    w_m \frac{\partial \phi^m}{\partial t} + \frac{\partial w}{\partial t}
    -
    w\frac{\partial \phi^m}{\partial t}w_m
    \bigg )
           \nonumber
    \\
    &= \frac{1}{P}w_t.
\end{align}
\end{proof}

\section{(Normalized) $\sigma_k$ curvature flow}
\label{sect:flowk}
In \cite{MR2647128}, Cabezas-Rivas and Sinestrari studied the volume preserving flow 
\begin{align}
    \label{flow_volumepreserve}
    X_t = (-\sigma_k(L)^{\alpha} + h(t))\nu,
\end{align}
where $\nu$ is the outer unit normal, $\alpha>\frac{1}{k}$, and $h(t)$ is the normalization factor defined as 
\begin{align}
    h(t) &:= \frac{1}{\text{Area}(M(t))}\int_{M(t)}\sigma_k(L)^{\alpha} d\mu_t.
\end{align}
This choice of $h$ ensures the volume remains fixed under the flow since
\begin{align}
    \frac{d}{dt}\text{Vol}(\Omega (t))
    &= \int_{M(t)} <X_t, \nu> d\mu_t
    \nonumber
    \\
    \label{volume_preserve_meancurv}
    &=\int_{M(t)}-\sigma_k(L)^{\alpha} + h(t) d\mu_t = 0. 
\end{align}
They assumed a pinching condition on $M(0)$ and showed that the condition is preserved along the flow. The pinching condition is a stronger condition than convexity, and their main theorem on the existence of a solution to the flow is summarized in the following theorem.

\begin{Thm}(Cabezas-Rivas and Sinestrari,  \cite{MR2647128}[Theorem 1.1])
\label{sinestrari_thm}
Given a fixed $\alpha>0$ and $1 \leq k \leq n$, there exists a constant $C_p \in (0,\frac{1}{n^n})$ depending only on $k, \alpha,$ and $n$, such that if the initial surface $M(0)$ has the property that at every point
\begin{align}
\label{pinching_condition}
    \sigma_n(L) > C_p \sigma_1^n(L)>0,
\end{align}
then there exists a unique solution of surfaces $\{ M(t)\}$ such that
\begin{enumerate}
    \item The pinching condition \eqref{pinching_condition} holds for all $t>0$ for which a solution to the flow exists,
    \item $M(t)$ exists for all $t\in (0, \infty)$,
    \item $\{ M(t)\} $ converges exponentially fast to a sphere and the volume is preserved along the flow.
\end{enumerate}
\end{Thm}

Note that the Newton-Maclaurin inequality $\frac{1}{n^n}\sigma^n_1(L)\geq\sigma_n(L)$ gives an opposite inequality of the pinching condition and attains equality when the surface is a sphere.
%but we get equality when the surface is a sphere.  
Thus with the established $C_p$, we know that a nearly spherical surface, with $\|u\| < \epsilon$ for a sufficiently small $\epsilon$, will satisfy the pinching condition. 
Additionally, any surface satisfying the pinching condition will be convex. 
In \cite{MR3830246}, Bertini and Sinestrari were able to remove the pinching condition when $k=2$, and they only assumed convexity on the initial surface to show the existence of a solution.

Just as with the inverse curvature flow, %we are concerned with the derivatives of various quantities along the flow to establish a quantitative quermassintegral inequality. We note that all the computations in 
using Proposition \ref{flow_derivatives_prop} for $\frac{1}{f}= -(\sigma_k^{\alpha}(L)- h)\nu$, one can compute 
$\frac{d}{dt}\int_{M(t)}\sigma_{k-1}(L)$.
%which is contained in the next lemma (See \cite{MR3830246}).

%which states that $\int_{M(t)}\sigma_{k-1}(L) d\mu_t$ is decreasing along the flow. 

\begin{Lem}(Bertini and Sinestrari, \cite{MR3830246}[Lemma 3.3])
\label{Bertini_deriv}
For a solution of hypersurfaces $\{M(t)\} $ along the flow \eqref{flow_volumepreserve},
\begin{align}
   & \frac{d}{dt}\int_{M(t)}\sigma_{k-1}(L)  d\mu_t= 
     - k\int_{M(t)}
    \sigma_{k}(\sigma_k^{\alpha}(L) - h )
    d\mu_t
\nonumber\\    
   =& - k\int_{M(t)}
    (\sigma_{k}(L) - h^{1/\alpha})(\sigma_k^{\alpha}(L) - h )
    d\mu_t
    \leq 0,
\end{align}
and there is equality only when $M(t)$ is a sphere.
\end{Lem}
Note that this lemma, along with \eqref{volume_preserve_meancurv}, proves the $(k-1,-1)$-quermassintegral inequality in the case where $M$ satisfies the pinching condition. We present the proof here for readers' convenience.
%\edz{It should be $(k-1,-1)$-quermassintegral inequality.See sidemark on Page 6.}
\begin{proof}
We find from the results of Proposition \ref{flow_derivatives_prop} that, setting $G:=
-\sigma_k^{\alpha}(L)+ h$,

%\edz{$\Sigma_{m-1,1}(h^i_j,(h^2)^i_j)$ not defined. Use $T_{m-1}$ notation instead. Be consistent with notation in Prop 3.2 after making the notation change.}
\begin{align}
    \frac{d}{dt}
    \int_{M(t)} \sigma_{m}(L) 
    d\mu_t
    &=
       \int_{M(t)} \partial_t(\sigma_{m}(L))
    d\mu_t
    +
       \int_{M(t)} \sigma_{m}(L) 
   \partial_t d\mu_t
   \nonumber
   \\
     &=
       -\int_{M(t)} 
       \nabla_j([T_{m-1}]^i_j(L) \nabla_i G)
       +G
       [T_{m-1}]^{i}_j(L) (L^2)^j_i
    d\mu_t 
    \nonumber
   \\ 
      &+\int_{M(t)} \sigma_{m}(L) 
   \sigma_1(L)Gd\mu_t.
\end{align}
%Then for the same reasoning as in 
Using $
[T_{m-1}]^i_j(L) (L^2)^j_i = 
\sigma_1(L) \sigma_{k-1}(L)- k \sigma_k(L),
$ we have
\begin{align}
    \frac{d}{dt}
    \int_{M(t)} \sigma_{k-1}(L) d\mu_t
    &=
   - k\int_{M(t)}
    \sigma_{k}(\sigma_k^{\alpha}(L) - h )
    d\mu_t.
    \end{align}
Observing that $-h^{1/\alpha}\int_{M(t)}
(\sigma_k^{\alpha}(L) - h )
    d\mu_t =0$, we find
    \begin{align}
      \frac{d}{dt}
    \int_{M(t)} \sigma_{k-1}(L) d\mu_t
    &=
   - k\int_{M(t)}
    (\sigma_{k}(L) - h^{1/\alpha})(\sigma_k^{\alpha}(L) - h )
    d\mu_t.
    \end{align}
   Since $(\sigma_{k}(L) - h^{1/\alpha})$ and $(\sigma_k^{\alpha}(L) - h )$ share the same sign, we find
    $$  \frac{d}{dt}
    \int_{M(t)} \sigma_{k-1}(L) d\mu_t \leq 0.$$
\end{proof}

\section{Stability along the inverse curvature flow}
In this section, we analyze the stability of the $(k,k-1)-$ quermassintegral inequalities for $k-$convex starshaped surfaces $M(t)$, along the flow 
\eqref{statement_inverse_flow}: %when $f = \frac{\sigma_k(L)}{\sigma_{k-1}(L)}$ . That is the flow
\begin{align}
\label{inverseflowstatement}
    X_t &= \frac{\sigma_{k-1}(L)}{\sigma_{k}(L)}\nu.
\end{align}
Writing $M(t)$ in spherical coordinates where $w$ is the radial function, $f = \frac{\sigma_k(L)}{\sigma_{k-1}(L)}$, we saw by Lemma \ref{deriv_of_radius} that
\begin{align}
\label{radius_flow_formula}
    w_t &= \frac{\sqrt{1 + \frac{|\nabla w|^2}{w^2}}}{f} .
\end{align}
We rescale the surfaces to $\Tilde{M}(t)$ so that the new radial function satisfies
 \begin{align}
 \label{scaledsurface}
     \Tilde{w} &= e^{-rt}w,
   \hspace{.15in}
   r = \frac{{n \choose k-1}}{{n \choose k}}.
 \end{align}
From \eqref{deriv_k} and \eqref{deriv_kminus1}, we saw that this choice of $r$ ensures $\int_{\Tilde{M}(t)}\sigma_{k}(\tilde{L}) d\mu_t$ is decreasing  and $\int_{\Tilde{M}(t)}\sigma_{k-1}(\tilde{L}) d\mu_t$ is constant along the flow.  We set
\begin{align}
    u &:= \Tilde{w}-1.
\end{align}

For the computations below, we take $M(t)$ to be nearly spherical so that $||u||_{C^2}< \epsilon$. Using the formula in Lemma 3.1 and Lemma 4.1 
in \cite{MR2441221}, %\ref{sigmak} and Lemma \ref{expansion},
and denoting $O(W):=
 O(u\Delta u)
    + O(u^2) + O(|\nabla u|^2) + O(|D^2u|^2)$, we find
\begin{align}
    \sigma_{k-1}(\Tilde{L})
    &=
    \frac{{n \choose k-1}}{((1+u)^2 + |\nabla u|^2)^{\frac{k+1}{2}}}
    \bigg (
    1+2u
    -\frac{k-1}{n}(1+u)\Delta u
    \bigg )
    + O(W)
    \nonumber
    \\
      &=
  {n \choose k-1}(1-(k+1)u)
    \bigg (
    1+2u
    -\frac{k-1}{n}\Delta u
    \bigg )
      + O(W)
    \nonumber
      \\
      \label{kminus1expan}
      &=
  {n \choose k-1}
    \bigg (
    1 - (k-1)u
    -\frac{k-1}{n}\Delta u
    \bigg )
    +O(W).
\end{align}
And,
\begin{align}
    \frac{1}{\sigma_{k}(\Tilde{L})}
    &= \frac{1}{{n \choose k}}\frac{1}{1 - (1 - \sigma_{k}(\Tilde{L})/{n \choose k})}
\nonumber
\\
      &=
      \frac{1}{{n \choose k}}\bigg 
      (1 + 
      \bigg(1 - \frac{\sigma_{k}(\Tilde{L})}{{n \choose k}} \bigg)
      + \sum_{j = 2}^{\infty}
      \bigg(1 - \frac{\sigma_{k}(\Tilde{L})}{{n \choose k}}
      \bigg)^j \bigg)
\nonumber    
\\
      &=
      \frac{1}{{n \choose k}}\bigg (2 - \frac{\sigma_{k}(\Tilde{L})}{{n \choose k}} \bigg) + O(W)
   \nonumber
   \\
       &=
      \frac{1}{{n \choose k}}\bigg (2 -
      (1 - (k+2)u)(1 + 2u - \frac{k}{n}\Delta u) \bigg) + 
      O(W)
\nonumber      \\
    \label{kinverseexpan}
      &=
      \frac{1}{{n \choose k}}\bigg (1 +ku + \frac{k}{n}\Delta u) \bigg)   
      + O(W).
\end{align}
These formulas are utilized in the following lemma. 
%\edz{Starting point of the 2nd time review.}
\begin{Lem}
\label{unormandak}
Suppose $\{\Tilde{M}(t)\}$ are the rescaled surfaces, defined as in \eqref{scaledsurface}, of a solution to the flow \eqref{inverseflowstatement}. If at some $t_0 \geq 0$,  $\Tilde{M}(t_0)$ is a nearly spherical surface where $||u(t_0)||_{C^2}<\epsilon$ for sufficiently small $\epsilon>0$, it holds that at $t_0$
\begin{align}
    \frac{d}{dt}||u||^2_{L^2}&= 
    \frac{-2}{n} \frac{{n \choose k-1}}{{n \choose k}}||\nabla u||^2_{L^2}
        +O(\epsilon)||u||_{W^{2,2}}^2, 
    \end{align}
    and
    \begin{align}
     \frac{d}{dt}||\nabla u||^2_{L^2}&=
     \frac{-2}{n} \frac{{n \choose k-1}}{{n \choose k}}||\Delta u||^2_{L^2}
         +O(\epsilon)||u||_{W^{2,2}}^2
     .
\end{align}
\end{Lem}
% \edz{1. Does $L^2$ of $u$ use $dA$, $\mu_t$ or $\tilde{\mu_t}$? 2. Does the gradient and $\Delta$ with respect to the metric of standard sphere? }
\begin{proof}

Combining the expansions of $\sigma_{k-1}(\Tilde{L})$  and $   \frac{1}{\sigma_{k}(\Tilde{L})} $ from \eqref{kminus1expan} and \eqref{kinverseexpan},
\begin{align}
   \frac{\sigma_{k-1}(\Tilde{L})}{\sigma_{k}(\Tilde{L})} 
      &=
 \frac{{n \choose k-1}}{{n \choose k}} 
    \bigg (
    1 - (k-1)u
    -\frac{k-1}{n}\Delta u
    \bigg )
   \bigg  (1 +ku + \frac{k}{n}\Delta u \bigg)
    +O(W)
      \nonumber
    \\
         &=
 \frac{{n \choose k-1}}{{n \choose k}} 
    \bigg (
    1 +u
    +\frac{1}{n}\Delta u
    \bigg )
    + O(W).
\end{align}
We use this expansion to compute 
\begin{align}
     \frac{d}{dt}||u||_{L^2}^2
    &= 2\intsphere uu_t \textit{dA}
\nonumber    \\
&= 2\intsphere u\tilde{w}_t \textit{dA}
\nonumber    \\
&= 2\intsphere u 
\bigg(-re^{-rt}w + e^{-rt}w_t \bigg) \textit{dA}
\nonumber    \\
&= 2\intsphere u \bigg(-\frac{{n \choose k-1}}{{n \choose k}}(1+u) 
+ e^{-rt} \sqrt{1 + \frac{|\nabla u|^2}{(1+u)^2}}\frac{\sigma_{k-1}({L})}{\sigma_{k}({L})} \bigg) \textit{dA}.
\end{align}
After rescaling, we have
%\edz{(5.11) Should be $O(\epsilon)||u||^2_{W^{2,2}} $}
\begin{align}
 \frac{d}{dt}||u||_{L^2}^2&= 2\intsphere u
\bigg(-\frac{{n \choose k-1}}{{n \choose k}}(1+u)  
+  \sqrt{1 + \frac{|\nabla u|^2}{(1+u)^2}}\frac{\sigma_{k-1}(\Tilde{L})}{\sigma_{k}(\Tilde{L})}
\bigg) \textit{dA}
   \nonumber
   \\
     &=
     2 \frac{{n \choose k-1}}{{n \choose k}}
     \intsphere u
     \bigg(
     -(1+u) 
   + 
    ( 1 +u
    +\frac{1}{n}\Delta u)
    \bigg)
    \textit{dA}
   +O(\epsilon)||u||^2_{W^{2,2}}
   \nonumber
     \\
     &=
       \frac{2}{n} \frac{{n \choose k-1}}{{n \choose k}}
     \intsphere u
     \Delta u
    \textit{dA}
   +O(\epsilon)||u||^2_{W^{2,2}}
   \nonumber
     \\
     &=
       \frac{-2}{n} \frac{{n \choose k-1}}{{n \choose k}}
     \intsphere |\nabla u|^2
    \textit{dA}
   +O(\epsilon)||u||^2_{W^{2,2}}.
\end{align}
%\edz{ Does the second equation in (5.11) use expansion formula of Lemma 4.1 in \cite{MR2441221}?}
Next, we find
\begin{align}
    \intsphere
     \frac{d}{dt} |\nabla w|^2 \textit{dA}
     &=
      \intsphere
     2w^j(w_t)_j
     \textit{dA}
    \nonumber   \\
      &=
     - \intsphere
     2w_t\Delta w
     \textit{dA}
     \nonumber   \\
      &=
    -  \intsphere
   2
     \Delta w\sqrt{1 + \frac{|\nabla u|^2}{(1+u)^2}} \frac{\sigma_{k-1}(\Tilde{L})}{\sigma_{k}(\Tilde{L})}
     \textit{dA}.
\end{align}
%\edz{write the last line in (5.12) as $|\nabla u|^2/u^2$; $G$ as $1/f$ to be consistent with (5.11). }
Therefore,
\begin{align}
     \intsphere
     \frac{d}{dt} |\nabla u|^2 \textit{dA}
     &=
      \intsphere
     \frac{d}{dt} \bigg( e^{-2rt}|\nabla w|^2 \bigg)\textit{dA}
     \nonumber \\
     &= \intsphere
     -2re^{-2rt}|\nabla w|^2
     +
     e^{-2rt} \frac{d}{dt} |\nabla w|^2 \textit{dA}
      \nonumber \\
     &= \intsphere
     -2\frac{{n \choose k-1 }}{{n \choose k}}|\nabla u|^2
     +
     e^{-2rt} \frac{d}{dt} |\nabla w|^2 \textit{dA}
         \nonumber \\
     &= \intsphere
     -2\frac{{n \choose k-1 }}{{n \choose k}}|\nabla u|^2
     -2
     e^{-2rt}\frac{\sigma_{k-1}(L)}{\sigma_k(L)}
     \bigg (
     \Delta w\sqrt{1 + \frac{|\nabla w|^2}{w^2}}
     \bigg)
     \textit{dA}.
 \end{align}
After rescaling, we find
 \begin{align}
    \intsphere
     \frac{d}{dt} |\nabla u|^2 \textit{dA}  &= \intsphere
     -2\frac{{n \choose k-1 }}{{n \choose k}}|\nabla u|^2
    -2
    \frac{\sigma_{k-1}(\tilde{L})}{\sigma_k(\tilde{L})}
     \bigg (
     \Delta u \sqrt{1 + \frac{|\nabla u|^2}{(1+u)^2}}
     \bigg)
     \textit{dA}
       \nonumber \\
     &=-2\frac{{n \choose k-1 }}{{n \choose k}} \intsphere
     |\nabla u|^2
    +
   \Delta u \bigg( 
    1+ u + \frac{1}{n}\Delta u
    \bigg)
        \textit{dA}
        + O(\epsilon)||u||^2_{W^{2,2}}    
        \nonumber \\
     &=-2\frac{{n \choose k-1 }}{{n \choose k}} \intsphere
     |\nabla u|^2
  +u\Delta u
  +
  \frac{1}{n} (\Delta u)^2
        \textit{dA}
        + O(\epsilon)||u||^2_{W^{2,2}}
         \nonumber \\
     &=-2\frac{{n \choose k-1 }}{{n \choose k}} \intsphere
  \frac{1}{n} (\Delta u)^2
        \textit{dA}
        + O(\epsilon)||u||^2_{W^{2,2}}.
\end{align}
%\edz{1. In the second equation in (5.14) we need to use Lemma 4.1 of \cite{MR2441221}.\\ 2. Should be $O(\epsilon)||u||^2_{W^{2,2}}$.}

\end{proof}
Next, we expand the derivative of$\int_{\Tilde{M}} \sigma_k(\Tilde{L}) d\mu_t$ along the flow. 

%\edz{ $\textit{dA}$ Should be $d\mu_t$. I've replaced all in Lemma 5.2.}
\begin{Lem}
\label{meanflow_sigma_deriv}
For the nearly spherical surface $\Tilde{M}(t_0)$ in Lemma \ref{unormandak}, at $t_0$ it holds that
\begin{align}
\label{first_sigmak_flow_ineq}
     \frac{d}{dt} \int_{\Tilde{M}} \sigma_k(\Tilde{L}) d\mu_t
      &\leq 
    (k+1) \frac{{n  \choose k-1}{n  \choose k+1}}{{n \choose k}} \bigg ( \frac{1}{n}||\nabla u||_{L^2}^2
    - \frac{1}{n^2}|| \Delta u||_{L^2}^2
    \bigg)
     +O(\epsilon)||u||_{W^{2,2}}^2.
\end{align}

\end{Lem}
\begin{proof}
From formula \eqref{deriv_k}, we find
\begin{align}
     &\frac{d}{dt} \int_{\Tilde{M}} \sigma_k(\Tilde{L}) d\mu_t\nonumber\\
      =&
     (k+1)\intsphere \frac{1}{\sigma_k(\Tilde{L})}
     \bigg ( 
     \sigma_{k-1}(\Tilde{L})\sigma_{k+1}(\Tilde{L})
     -
     \frac{{n \choose k-1}{n \choose k+1}}{{n \choose k}^2}\sigma_k^2(\Tilde{L})
     \bigg) \sqrt{ \text{det } \tilde{g}}
     \textit{dA},
\end{align}
and we found in \eqref{area_element} that
\begin{align}
     \sqrt{ \text{det } \tilde{g}} &=(1+u)^n\sqrt{1 + \frac{|\nabla u|^2}{(1+u)^2}}.
\end{align}
From the Newton-Maclaurin inequality 
 $ \sigma_{k-1}(\Tilde{L})\sigma_{k+1}(\Tilde{L})
     -
     \frac{{n \choose k-1}{n \choose k+1}}{{n \choose k}^2}\sigma_k^2(\Tilde{L}) \leq 0$, formula \eqref{kinverseexpan} and Lemma 4.1 of \cite{MR2441221}, we conclude 
\begin{align}
     \frac{d}{dt} \int_{\Tilde{M}} \sigma_k(\Tilde{L}) d\mu_t
      &\leq 
     \frac{k+1}{{n \choose k}}\intsphere 
     \sigma_{k-1}(\Tilde{L})\sigma_{k+1}(\Tilde{L})
     -
     \frac{{n \choose k-1}{n \choose k+1}}{{n \choose k}^2}\sigma_k^2(\Tilde{L})
     \textit{dA}
  \nonumber
  \\  
  &+ 
     O(\epsilon)\intsphere 
     \sigma_{k-1}(\Tilde{L})\sigma_{k+1}(\Tilde{L})
     -
     \frac{{n \choose k-1}{n \choose k+1}}{{n \choose k}^2}\sigma_k^2(\Tilde{L})
     \textit{dA}.
\end{align}
Now we expand the expression  $   \sigma_{k-1}(\Tilde{L})\sigma_{k+1}(\Tilde{L})
     -
     \frac{{n \choose k-1}{n \choose k+1}}{{n \choose k}^2}\sigma_k^2(\Tilde{L})$ from the formula in Lemma 3.1 in \cite{MR2441221}. %\ref{sigmak}. 
     We find
     $
      \sigma_{k-1}(\Tilde{L})\sigma_{k+1}(\Tilde{L})
((1+u)^2 + |\nabla u|^2)^{k+2}
$
is equal to the product:
    \begin{align}
  & 
     \bigg (
     \sum_{m=0}^{k-1}\frac{(-1)^m{n -m\choose k-1-m}}{(1+u)^m}
     \bigg(
     (1+u)^2\sigma_m(D^2u) + \frac{n+k-1-2m}{n-m} u^iu_j[T_m]^j_i(D^2u)
     \bigg)
     \bigg)
     \nonumber
     \\
     \times &
     \bigg (
      \sum_{m=0}^{k+1}\frac{(-1)^m{n -m\choose k+1-m}}{(1+u)^m}
     \bigg(
     (1+u)^2\sigma_m(D^2u) + \frac{n+k+1-2m}{n-m} u^iu_j[T_m]^j_i(D^2u)
     \bigg)
     \bigg),
      \end{align} 
     and $ \sigma_k^2(\Tilde{L})
((1+u)^2 + |\nabla u|^2)^{k+2}$ equals:
      \begin{align}
  &   \bigg(
    \sum_{m=0}^{k}\frac{(-1)^m{n -m\choose k-m}}{(1+u)^m}
     \bigg(
     (1+u)^2\sigma_m(D^2u) + \frac{n+k-2m}{n-m} u^iu_j[T_m]^j_i(D^2u)
     \bigg)
     \bigg)^2.
    \end{align} 
 We collect the coefficients occurring in front of the lower order terms in the expression
     \begin{align}
       ( (1+u)^2 + |\nabla u|^2)^{k+2}(\sigma_{k-1}(\Tilde{L})
      \sigma_{k+1}(\Tilde{L})
     -
     \frac{{n \choose k-1}{n \choose k+1}}{{n \choose k}^2}\sigma_k^2(\Tilde{L})),
     \end{align}
for which we have the following:
    \begin{itemize}
        \item $(1+u)^4: 0$
         \item $(1+u)^2|\nabla u|^2: 0$
         \item $(1+u)^3 \Delta u :0$
           \item $(1+u)^2 \sigma_2(D^2u) : {n  \choose k-1}{n  \choose k+1}\frac{2}{n(n-1)}$
           \item $(1+u)^2(\Delta u)^2:  -{n  \choose k-1}{n  \choose k+1}\frac{1}{n^2}$
    \end{itemize}
   The rest of the terms are in $O(\epsilon)u^2 +O(\epsilon)|\nabla u|^2+O(\epsilon)|D^2u|^2 $. Because $$\frac{1}{((1+u)^2 + |\nabla u|^2)^{k+2}}= 1 + O(\epsilon),$$ we obtain $ \sigma_{k-1}(\Tilde{L})\sigma_{k+1}(\Tilde{L})
     -
     \frac{{n \choose k-1}{n \choose k+1}}{{n \choose k}^2}\sigma_k^2(\Tilde{L})$ is equal to:
    \begin{align}
     & {n  \choose k-1}{n  \choose k+1}
      \bigg ( \frac{2}{n(n-1)}\sigma_2(D^2u)
    -
      \frac{1}{n^2}(\Delta u)^2 \bigg)
     +O(\epsilon)u^2 +O(\epsilon)|\nabla u|^2
     \nonumber\\
   &  +O(\epsilon)|D^2u|^2.
    \end{align}
    Hence,
\begin{align}
     \frac{d}{dt} \int_{\Tilde{M}} \sigma_k(\Tilde{L}) d\mu_t
      \leq &
    (k+1) \frac{{n  \choose k-1}{n  \choose k+1}}{{n \choose k}}\intsphere 
    \frac{2}{n(n-1)}\sigma_2(D^2u)
    -
    \frac{1}{n^2}(\Delta u)^2
     \textit{dA}\nonumber \\
    & +O(\epsilon)||u||_{W^{2,2}}^2 .
\end{align}

In the proof of Lemma 4.2 of \cite{MR2441221},
%\ref{expansion}, 
we used integration by parts to find $\intsphere 
    \sigma_2(D^2u)  \textit{dA}=
    \frac{n-1}{2}\intsphere 
    |\nabla u|^2  \textit{dA} +O(\epsilon)||\nabla u||_{L^2}^2$. Hence
\begin{align}
     \frac{d}{dt} \int_{\Tilde{M}} \sigma_k(\Tilde{L}) d\mu_t
      &\leq 
         \frac{(k+1){n  \choose k-1}{n  \choose k+1}}{n\cdot{n \choose k}} \bigg ( ||\nabla u||_{L^2}^2
    - \frac{1}{n}|| \Delta u||_{L^2}^2
    \bigg)
     +O(\epsilon)||u||_{W^{2,2}}^2.
\end{align}

\end{proof}

In Lemma 4.2 of \cite{MR2441221},
%\ref{spherharmbound},
we obtained a Poincar\'{e} inequality when the barycenter of a surface is at the origin. Following a similar argument to this lemma, we now obtain a similar inequality, but instead we compare $||\Delta u||_{L^2}^2$ and $||\nabla u||_{L^2}^2$. We also relax the conditions on the barycenter slightly.

\begin{Lem}
\label{poincare_laplace_bound}
Suppose $M$ is a hypersurface in $\mathbb{R}^{n+1}$ which is starshaped with respect to the origin, so that $M = \{(1+u(x))x: x\in \partial B\}$ where $u \in C^2(\partial B)$. Further assume that $M$ is nearly spherical with $||u||_{C^2} < \epsilon$, and for some fixed $K>0$ the barycenter of $M$ satisfies $|\text{bar}(M)|^2 \leq K\epsilon ||u||_{W^{2,2}}^2$. Then, 
\begin{align}
     ||\Delta u||_{L^2}^2
     & \geq   2(n+1)
   ||\nabla u||_{L^2}^2
    -
  K'\epsilon||u||_{W^{2,2}}^2,
\end{align}
where $K'>0$ depends on the choice of $K$.
\end{Lem}
\begin{proof}
Write  $u = \sum_{k=0}^{\infty}a_kY_k$. From the proof of Lemma 5.2 in \cite{MR2441221},
%\ref{spherharmbound}, 
observe that if $|\text{bar}(M)|^2 \leq K\epsilon ||u||_{W^{2,2}}^2$, then for some $K'>0$ (depending on $K$) we have $a_1^2 \leq K'\epsilon||u||_{W^{2,2}}^2$. Furthermore, we find that
\begin{align}
   ||\Delta u||_{L^2}^2&= \sum_{k=1}^{\infty}\lambda_k^2a_k^2
\nonumber
\\
 &\geq
 |\lambda_2|
    \sum_{k=1}^{\infty}
    |\lambda_k|a_k^2
    +
    (|\lambda_1|-|\lambda_2|)|\lambda_1|a_1^2
    \nonumber
    \\
    \label{poincare_laplacian}
    &=
    2(n+1)
   ||\nabla u||_{L^2}^2
      +
    (|\lambda_1|-|\lambda_2|)|\lambda_1|a_1^2.
\end{align}
\end{proof}

We recall the following proposition from \cite{MR2441221}.

\begin{Prop}(\cite{MR2441221}[Proposition 5.1])
\label{propgen}
 Fix  $j$ where $ 0 \leq j < k$. Suppose $\Omega = \{ (1+u(\frac{x}{|x|}))x: x \in B \} \subseteq \mathbb{R}^{n+1}$, where $u \in C^3(\partial B)$, $I_j(\Omega) = I_j(B)$, and $\text{bar}(\Omega) = 0$.
Assume for sufficiently small $\epsilon >0$ that $||u||_{W^{2,\infty}}<\epsilon$. Then,
\begin{align} I_k(\Omega) - I_k(B) \geq&
 {n  \choose k}\frac{(n-k)(k-j)}{2n}
 \bigg( 
 \bigg(1+ O(\epsilon) \bigg) || u||_{L^2}^2 
 \nonumber\\
 &+ \bigg(\frac{1}{2} +O(\epsilon)
 \bigg )||\nabla u||_{L^2}^2 
 \bigg ).
   \end{align} 
\end{Prop}

Now we are fully equipped with the formulas needed to prove the next proposition. Under the conditions of Proposition \ref{propgen}, when $I_{k-1}(\Omega)) = I_{k-1}(B)$ we have 
\begin{align} 
\label{gen_prop_rephrased}I_k(\Omega) - I_k(B) \geq
(1 + O(\epsilon))A(t),
   \end{align} 
   where we recall the notation of \eqref{A}
 \begin{align}
     A(t) &:= A(\Omega(t))=
      {n  \choose k}\frac{(n-k)}{2n}
 \bigg( 
 || u||_{L^2}^2  + \frac{1}{2}||\nabla u||_{L^2}^2 
 \bigg ) .
 \end{align}
%The quantitative quermassintegral inequalities that we established for nearly spherical sets quickly followed from this proposition. 
The rest of this section is devoted to providing a new approach to proving the quantitative quermassintegral inequalities \eqref{gen_prop_rephrased} by comparing $\frac{d}{dt}(I_k(\tilde{\Omega})(t) - I_k(B)) $ and $\frac{d}{dt}
    A(t)$ along the flow.

For convenience, we recall the statement of Theorem \ref{main_inversemeancurv_inequality}.

\begin{Thm}
Suppose $M(t)$ is a solution to the flow \eqref{inverseflowstatement}) and
$\Tilde{M}(t)$ is the rescaled surface in \eqref{scaledsurface} of $M(t)$. 
Additionally, assume at $t_0$ that $M(t_0)$ is nearly spherical with $||u(t_0)||_{W^{2,\infty}} < \epsilon$ and that the barycenter of $\Tilde{M}(t_0)$ satisfies $|bar(\tilde{M}(t_0))| \leq K \epsilon ||u(t_0)||_{W^{2,2}}^2$ for fixed a $K>0$. Then, for any small $\eta>0$,
\begin{align}
    \frac{d}{dt}(I_k(\tilde{\Omega}(t_0)) - I_k(B))
    \leq 
    (1-\eta)\frac{d}{dt}
    A(t_0),
\end{align}
and the choice of a sufficiently small $\epsilon>0$ depends on $\eta$ and $K$.

Moreover, along any solution to the flow \eqref{inverseflowstatement} where
$|bar(\tilde{M}(t))| \leq K \epsilon ||u||_{W^{2,2}}^2$
\textit{holds for sufficiently large $t$, we have}
\begin{align}
   \liminf_{t\rightarrow \infty}  \frac{I_k(\tilde{\Omega}(t)) - I_k(B)}
     {A(t)}
   \geq 1 .
\end{align}
\end{Thm}

\begin{Rem}
If the initial surface $M(0)$ is n-symmetric (symmetric with respect to reflection over each coordinate axis), then $M(t)$ remains n-symmetric throughout the flow, in which case the barycenter remains at the origin during the entire flow, thus satisfying the conditions on the barycenter for this theorem. %One major defect of our result is that we cannot prove the conditions on the barycenter are satisfied along the flow if it was initially satisfied on $M(0)$. 

\end{Rem}

\begin{proof}
Lemma \ref{unormandak} immediately yields
\begin{align}
    \frac{d}{dt}A(t)
    &=
      - {n  \choose k}\frac{n-k}{2n}
       \frac{2}{n}\frac{{n \choose k-1}}{{n \choose k}}
 \bigg( 
 || \nabla u||_{L^2}^2  + \frac{1}{2}||\Delta u||_{L^2}^2 
 \bigg )
 +O(\epsilon)
 ||u||_{W^{2,2}}^2
 \nonumber
 \\
 &=
     -  \frac{(k+1)}{n^2}\frac{{n \choose k-1} {n  \choose k+1}}{{n \choose k}}
 \bigg( 
 || \nabla u||_{L^2}^2  + \frac{1}{2}||\Delta u||_{L^2}^2 
 \bigg )
 +O(\epsilon)
 ||u||_{W^{2,2}}^2.
\end{align}
Combining the inequalities in \eqref{first_sigmak_flow_ineq} and \eqref{poincare_laplacian} yields
\begin{align}
     \frac{d}{dt} \int_{\Tilde{M}} \sigma_k(\Tilde{L}) d\mu_t
      \leq &
    \frac{(k+1)}{n^2}
    \frac{{n  \choose k-1}{n  \choose k+1}}{{n \choose k}}
  \bigg(  
    n||\nabla u|_{L^2}^2 - \frac{1}{2}||\Delta u||_{L^2}^2
   - \frac{1}{2}||\Delta u||_{L^2}^2
    \bigg) \nonumber\\
    &+O(\epsilon)||u||_{W^{2,2}}^2.
    \nonumber
    \\
    \label{inequality_before_eta}
     \leq &
    \frac{(k+1)}{n^2}
    \frac{{n  \choose k-1}{n  \choose k+1}}{{n \choose k}}
  \bigg( 
    -||\nabla u||_{L^2}^2
   - \frac{1}{2}||\Delta u||_{L^2}^2
    \bigg) +O(\epsilon)||u||_{W^{2,2}}^2.
\end{align}
At this point, we can almost conclude the first part of the theorem, but the term $O(\epsilon)||u||_{W^{2,2}}^2$ in both inequalities above does not allow for the conclusion that $\frac{d}{dt}(I_k(E) - I_k(B)) \leq \frac{d}{dt}A(t)$. However, if we multiply the left-hand side of \eqref{inequality_before_eta} by  $\frac{1}{1-\eta}$, we find
\begin{align}
\label{constant_biggerthanone_ineq}
      \frac{1}{1-\eta}\frac{d}{dt} \int_{\Tilde{M}} \sigma_k(\Tilde{L}) d\mu_t
     &\leq 
    \frac{(k+1)}{n^2}
    \frac{{n  \choose k-1}{n  \choose k+1}}{{n \choose k}}
  \bigg( 
    -||\nabla u|_{L^2}^2
   - \frac{1}{2}||\Delta u||_{L^2}^2
    \bigg)
    \nonumber
    \\
    &-C'||\Delta u||_{L^2}^2 +O(\epsilon)||u||_{W^{2,2}}^2,
\end{align}
where $C>0$ depends only on the choice of $\eta$.

Furthermore, noting that on the sphere $||\Delta u||_{L^2}^2=||D^2u||_{L^2}^2+(n-1)||\nabla u||_{L^2}^2$ and the inequalities in Lemma 5.2 of \cite{MR2441221} and Lemma \ref{poincare_laplace_bound}, %\ref{spherharmbound}
we find that $C'||\Delta u||_{L^2}^2$ dominates $O(\epsilon)||u||_{W^{2,2}}^2$. Thereby, for any choice of $0< \eta <1$, there is an $\epsilon>0$ where $||u||_{C^2}< \epsilon$ ensures
\begin{align}
\label{lemma_514_before_divide}
      \frac{d}{dt}(I_k(\tilde{\Omega}(t)) - I_k(B))
    \leq 
   (1-\eta) \frac{d}{dt}
    A(t).
\end{align}
Both sides of the inequality above are negative quantities, which gives
\begin{align}
     \frac{\frac{d}{dt}(I_k(\tilde{\Omega}(t)) - I_k(B))}
     {\frac{d}{dt}A(t)}
   \geq 1 - \eta.
\end{align}
Next, as a consequence of the exponential convergence of $M(t)$ to a sphere shown in \cite{MR1082861} and \cite{MR1064876}, we have that $||u||_{C^2}$ converges to $0$ as $t\rightarrow \infty$. Hence,
\begin{align}
   \liminf_{t\rightarrow \infty}  \frac{\frac{d}{dt}(I_k(\tilde{\Omega}(t)) - I_k(B))}
     {\frac{d}{dt}A(t)}
   \geq 1 .
\end{align}
Because both $I_k(M(t)) - I_k(B)$ and $A(t)$ approach $0$, we obtain, from a generalized version of L'Hopital's Rule, that 
\begin{align}
   \liminf_{t\rightarrow \infty}  \frac{I_k(\tilde{\Omega}(t)) - I_k(B)}
     {A(t)}
   \geq 1 .
\end{align}
\end{proof}

Before we can use the previous proposition to give a different proof of Proposition \ref{propgen}, at least in the case where the barycenter remains close enough to the origin, we need the following lemma. We need to show that if $\tilde{M}(0)$ is nearly spherical initially, then $\Tilde{M}(0)$ remains nearly spherical, which will allow us to apply the previous lemma for all $\Tilde{M}(t)$ along the flow. 

\begin{Lem}
\label{remain_nearly_spherical}
Along the flow \eqref{inverseflowstatement}, for any $\epsilon> 0$ there is a $\delta>0$ such that when $||u(0)||_{C^2}< \delta$, it holds that $||u(t)||_{C^2}<\epsilon$ for all $t\geq 0$. 
\end{Lem}

\begin{proof}
In this proof, we primarily examine some results from \cite{MR1082861} to help prove the lemma. In this paper, Urbas normalizes $f$ in the flow so that $f(\delta^{i}_j) =1$. This normalization would require that the rescaling to $\tilde{M}(t)$ satisfies $\tilde{w}(t) = e^{-t}w$ (instead of $e^{-rt}w$). In this proof, we continue using our conventions from this section and don't use their normalization. Note that in our rescaling we have $\tilde{f} = e^{rt} f$, where we recall $f = \frac{\sigma_k(L)}{\sigma_{k-1}(L)}$. We change constants from their proof to appropriately adapt to our notation. 
\cite{MR1082861}[Lemma 3.1] proves that
%\edz{Is this the $\omega$ defined in our flow or the flow in Urbas's paper?}
\begin{align}
   \min_{\partial B} w(0)\leq w(t)e^{-rt} \leq 
     \max_{\partial B} w(0),
\end{align}
which immediately yields
$  \min_{\partial B} u(0)\leq u(t)\leq  \max_{\partial B} u(0)$.
%\edz{Since $u=\tilde{\omega}-1$, should it be $|u(0)|\leq  \max_{\partial B} |u(0)|+ 2$?}
In \cite{MR1082861}[Lemma 3.2], Urbas showed
\begin{align}
    \frac{|\nabla w(t)|}{w(t)}
    &
    \leq 
 \max_{\partial B}  \frac{|\nabla w(0)|}{w(0)},
\end{align}
which yields
\begin{align}
    |\nabla u(t)| \leq   \frac {1+u(t)}{ \min_{\partial B}( 1+u(0))} 
    \max_{\partial B}|\nabla u(0)|
     \leq   \frac { \max_{\partial B} (1+u(0)) }{ \min_{\partial B}( 1+u(0))} 
    \max_{\partial B}|\nabla u(0)|
\end{align}
%\edz{Add in (5.41) the bound for numerator $1+ u(t)\leq \max_{\partial B} \omega(0) $}
At this point, we have concluded that $||u(t)||_{C^1} < \epsilon$ for some chosen $\delta>0$ on the initial conditions, and next we find bounds on the principal curvatures to obtain the desired bound on $||D^2u||_{L^{\infty}}$. At the end of \cite{MR1082861}[Lemma 3.3], Urbas showed
\begin{align}
\label{hmaxbound}
 h_{max}(t) &\leq h_{max}(0) -2rt,
\end{align} 
where  $h = \log \bigg( \frac{\kappa}{<X,\nu>} \bigg)$ and $h_{max}(t)$ is the maximum value taken over all points on the surface $M(t)$ and the principal curvatures $\kappa$. 
When rescaling the inputs, we observe
\begin{align}
    h &=\log \bigg( \frac{e^{-rt}\tilde{\kappa}}{<e^{rt}\tilde{X},\nu>} \bigg)
    \nonumber
    \\
     &=\log \bigg( e^{-2rt}\frac{\tilde{\kappa}}{<\tilde{X},\nu>} \bigg).
\end{align}
Because of the bound on $||u||_{C^1}$, we have that $<\tilde{X},\nu>$ can be made close to $1$ for small $\delta$. Then, given any $\beta>0$ we can ensure, for small enough $\delta$, that
\begin{align}
e^{-2rt}(\tilde{\kappa}_{max}(t) - \beta)
&\leq 
e^{-2rt}(\tilde{\kappa}_{max}(0) + \beta),
\end{align} 
where we have used that the inequality \eqref{hmaxbound} can be rewritten as $ e^{h_{max}(t)} \leq e^{-2rt}e^{h_{max}(0)}$. This inequality quickly simplifies to
\begin{align}
\tilde{\kappa}_{max}(t) 
&\leq 
\tilde{\kappa}_{max}(0) + 2\beta.
\end{align} 
For sufficiently small $||u(0)||_{C^2}$,  we find $\tilde{\kappa}_{max}(0) $ is close to 1 so that
\begin{align}
    \tilde{\kappa}_{max}(t) &\leq 1 + 3\beta.
\end{align}
Next we find a similar lower bound on the principal curvatures. In \cite{MR1082861}[Lemma 3.5], Urbas sets $G:=\frac{{n \choose k}}{{n \choose k-1}}\frac{\sqrt{1 + \frac{|\nabla w|^2}{w^2}}}{wf}$, and they conclude
\begin{align}
     %\mymathop{ min}
     \min_{\partial B} G (0)
     \leq G \leq  %\mymathop{
     \max_{\partial B} G(0).
\end{align}
Using $\tilde{f} = e^{rt}f$, we rewrite the above inequality as
\begin{align}
\frac{{n \choose k}}{{n \choose k-1}} \frac{ \sqrt{1 + \frac{|\nabla \tilde{w}|^2}{\tilde{w}^2}}}{\tilde{w}  } \frac{1}{\max_{\partial B} G (0)}
     \leq
     \tilde{f} \leq  
        \frac{{n \choose k}}{{n \choose k-1}}\frac{ \sqrt{1 + \frac{|\nabla \tilde{w}|^2}{\tilde{w}^2}}}{\tilde{w}  } 
     \frac{1}{\min_{\partial B} G(0)}.
\end{align}
%\edz{In (5.48), the denominator and numerator of $\frac{\tilde{w}}{\sqrt{1 + \frac{|\nabla \tilde{w}|^2}{\tilde{w}^2}}} $ are upside down.}
We find for sufficiently small $\delta$ that 
\begin{align}
 1- \beta
       \leq
    \frac{{n \choose k-1}}{{n \choose k}} \tilde{f} 
     \leq
      1+ \beta.
\end{align}
for all $t\geq 0$. Now, using the inequality  $ \frac{1}{n}\sigma_1(\tilde{L}) \geq \frac{{n \choose k-1}}{{n \choose k}}\tilde{f}$ (see, for example \cite{MR1465184}[Lemma 15.13]), we sum over the principal curvatures at any point on $M(t)$ to find
$\frac{1}{n}\sum_{i=1}^n \tilde{\kappa}_i
       \geq 1 - \beta$. Therefore,
   \begin{align}
    \tilde{\kappa}_1
      & \geq n  - \sum_{i=2}^n \tilde{\kappa}_i - n\beta
       \nonumber
      \\
       & \geq n  - (n-1)(1 + 3\beta) - n\beta = 1 - (4n-3)\beta.
   \end{align}
In conclusion, the principal curvatures remaining close to 1, together with $||u(t)||_{C^1}$. These ensure that $||u(t)||_{C^2} < \epsilon$.
\end{proof}

We are now ready to give a proof of the following corollary, which gives a new proof of Proposition \ref{propgen} when $j=k-1$ in the case that $M$ is n-symmetric (that is, $M$ is preserved under reflection over each coordinate axis). When setting $M(0) =M$, the symmetry condition will be preserved throughout the entire flow, thus keeping the barycenter of $M(t)$ at the origin for all $t\geq 0$. This allows us to apply Theorem \ref{main_inversemeancurv_inequality} to prove Corollary \ref{corollary_inversemean}.
\begin{comment}
\begin{Cor}
\label{corollary_inversemean}
Given any $\eta>0$, there is an $\epsilon>0$ such that any smooth $n$-symmetric, nearly spherical set $M$ that, where $M = \partial \Omega$, satisfies the inequality
\begin{align}
    I_k(\Omega) - I_k(B) \geq 
    (1- \eta)A, 
\end{align}
 when $||u||_{C^2}<\epsilon$ and $I_{k-1}(\Omega) = I_{k-1}(B)$.
\end{Cor}
\end{comment}

\begin{proof} of Corollary \ref{corollary_inversemean}.
Suppose $M$ is the initial surface, and $M(t)$ is the flow defined by \eqref{inverseflowstatement}. We denote the rescaled flow \eqref{scaledsurface} by $\tilde{M}(t) = \partial \tilde{\Omega}(t)$. Next, we set
\begin{align}
    S(t):= \frac{I_k(\tilde{\Omega}(t))- I_k(B)}{A(t)}.
\end{align}
For any $\eta>0$, we find $\epsilon>0$,  as in Theorem \ref{main_inversemeancurv_inequality}, so that $\frac{d}{dt}I_k(\tilde{\Omega}(t)) \leq (1-\eta)\frac{d}{dt}A(t)$ when $||u(t)||_{C^2} < \epsilon$. From Lemma \ref{remain_nearly_spherical}, we know for small enough $||u(0)||_{C^2}$ that $M(t)$ remains nearly spherical throughout the flow to ensure $\frac{d}{dt}I_k(\tilde{\Omega}(t)) \leq (1-\eta)\frac{d}{dt}A(t)$. We aim to prove the conclusion of this corollary for $M(0)$, which we note is the same as $\tilde{M}(0)$. 

We find
\begin{align}
    \frac{d}{dt}S(t)
    &=\frac{A(t)\frac{d}{dt}I_k(\tilde{\Omega}(t))
    -(I_k(\tilde{\Omega}(t))- I_k(B))\frac{d}{dt}A(t)
    }{A^2(t)}
    \nonumber
    \\
    &\leq 
  \frac{\frac{d}{dt}A(t)}{A^2(t)} ((1-\eta)\cdot A(t) 
    -(I_k\tilde{\Omega}(t))- I_k(B)).
\end{align}
Given any $t_0 \geq 0$, there are two cases:
\begin{enumerate}
    \item $S(t_0) \geq 1-\eta$. This is equivalent to $(1-\eta)\cdot  A(t_0) \leq I_k(\tilde{\Omega}(t_0))- I_k(B)$.
    \item $S(t_0) <1-\eta$, which implies $S$ is decreasing at $t_0$ because in this case $ (1-\eta)\cdot A(t_0)
    -(I_k(\tilde{\Omega}(t_0))- I_k(B)) > 0$, and we always have  $\frac{\frac{d}{dt}A(t_0)}{A^2(t_0)} \leq 0$ for small $||u||_{C^2}$.
\end{enumerate}

Suppose when $t_0=0$ we have case 1. This means $S(0) \geq 1-\eta$, which concludes this corollary for $M(0)$.

Now we consider what would happen if case 2 were to occur at $t_0=0$, for which we will find a contradiction and conclude that case 1 must occur at $t_0=0$. If case 2 were to occur, we would have $S(0) < 1-\eta$ and that $S$ begins decreasing initially. This implies that $S(t)< 1 -\eta$ for all $t>0$. So, $S(t)$ remains in case 2 for all $t>0$ and is thereby decreasing throughout the entire flow. This is a contradiction, however, since that implies $\liminf_{t\rightarrow \infty} S(t) \leq 1-\eta < 1$, but Theorem \ref{main_inversemeancurv_inequality} shows that $\liminf_{t\rightarrow \infty} S(t) \geq 1$. 

\end{proof}

\section{Stability along the (normalized) $\sigma_k$ curvature flow}

In this section, we examine stability of the $(k,-1)$-quermassintegral inequality of the flow of surfaces $M(t)$ along  \eqref{flow_volumepreserve} with $\alpha =1$, which is the flow
\begin{align}
\label{MeanCurvFlowStatement2}
    X_t &=  (-\sigma_k(L) + h(t))\nu,
\end{align}
where $h(t) = \frac{1}{\text{Area}(M(t))}\int_{M(t)}\sigma_k(L) \textit{d}\mu_t$.
The volume along the flow is preserved. Therefore, when studying stability for nearly spherical sets, the Fraenkel asymmetry can be approximated above by $||u||_{L^2}^2$, which we observed using formula (110) in \cite{MR2441221}.
%\eqref{holder_frankael}.

%\edz{Do we want to give a separate explanation of this formula (110) of \cite{MR2441221} in a lemma?}

So, next we differentiate the quantities $ ||u||_{L^2}^2$ and $ \frac{d}{dt}
     \int_{M(t)} \sigma_{k-1}(L) d\mu_t $ along the flow. 

\begin{Lem}
\label{MeanFlowDefEst}
For sufficiently small $\epsilon>0$, if $||u(t)||_{C^2} < \epsilon$ along the flow \eqref{MeanCurvFlowStatement2}, then 
\begin{align}
\label{uL2meancurv}
    \frac{d}{dt} ||u||_{L^2}^2
    &=  2 \frac{k}{n} {n \choose k} \intsphere
    nu^2 - |\nabla u|^2
    \textit{dA}
    +
    O(\epsilon)||u||_{W^{2,2}}^2,
\end{align}
and
\begin{align}
\label{sigmakmeancurv}
     \frac{d}{dt}
     \int_{M(t)} \sigma_{k-1}(L) \textit{d}\mu_t
     =&
       k^3{n \choose k}^2\intsphere
    -u^2 
    -
    \frac{1}{n^2}
    (\Delta u)^2
    +\frac{2}{n}
    |\nabla u|^2
    \textit{dA}\nonumber\\
    &+
    O(\epsilon)||u||_{W^{2,2}}^2.
\end{align}

\end{Lem}

\begin{proof}
We begin by estimating the function $G:= -\sigma_k(L) + h(t)$. 
To compute 
$$h(t): = \frac{1}{\text{Area}(M(t))}
\int_{M(t)} \sigma_k(L) d \mu_t,$$ 
recall from the proof of Proposition 4.3 of \cite{MR2441221}
%\ref{propvol} 
that when $V(\Omega(t)) = V(B)$, we have $ \int_{M(t)} \sigma_k(L) d\mu_t$ equals
\begin{align}
    &
{n \choose k}\text{Area}(\partial B) 
+{n \choose k} 
\frac{(n-k)(k+1)}{2n}
 \bigg( ||\nabla u||_{L^2}^2
  - n|||u||_{L^2}^2
  \bigg)
    + O(\epsilon)||u||_{W^{1,2}}^2.
\end{align}
%\edz{In (6.4), is coefficient in front of $\|u\|^2_{L^2}$ $-n$, not $n$? }
Additionally
%\begin{equation}
\begin{align}
    \text{Area}(M(t))
    =& \int_{\partial B}
    (1+u)^{n-1}\sqrt{(1+u)^2 + |\nabla u|^2}
     \textit{dA}
    \nonumber
    \\
    =&
    \int_{\partial B}
    (1+ (n-1)u+ \frac{(n-1)(n-2)}{2}u^2)(1 + u + \frac{1}{2}|\nabla u|^2)
    \textit{dA}\nonumber
      \\
    &+ O(\epsilon)||u||^2_{W^{1,2}}\nonumber
      \\
    =&
    \int_{\partial B}
    1 + nu + \frac{n(n-1)}{2}u^2+\frac{1}{2}|\nabla u|^2
    \textit{dA}
    + O(\epsilon)||u||^2_{W^{1,2}}.\label{6.6}
\end{align}
% \edz{In (6.5) $O(\epsilon)||u||_{W^{1,2}}$ should be  $O(\epsilon)||u||^2_{W^{1,2}}$.}
%\end{equation}
%Using that $\text{Vol}(\Omega)=\text{Vol}(\partial B)$
Using the assumption that $\text{Vol}(\Omega(t)) = \text{Vol}(B)$, we have 
%from formula \eqref{volume_formula} for the volume
that
\begin{align}
\label{volume_sub}
\intsphere  u \textit{ dA} = \intsphere \frac{-n}{2}u^2 \textit{ dA} +  O(\epsilon)||u||^2_{L^2}. 
\end{align}
We substitute \eqref{volume_sub} to \eqref{6.6}
to find
\begin{align}
     \text{Area}(M(t))
     &= \int_{\partial B}1 - \frac{n}{2}u^2+ \frac{1}{2}|\nabla u|^2 \textit{dA}+
      O(\epsilon)||u||^2_{W^{1,2}}.
\end{align}
 %\edz{In (6.7) $O(\epsilon)||u||_{W^{1,2}}$ should be $O(\epsilon)||u||^2_{W^{1,2}}$.}
Hence
\begin{align}
    \frac{1}{\text{Area}(M(t))}
    &=
    \frac{1}{\text{Area}(\partial B)}
    \frac{1}{1 - \frac{1}{\text{Area}(\partial B)}
    (\frac{n}{2}|||u||_{L^2}^2 - \frac{1}{2}||\nabla u||_{L^2}^2)
    +
    O(\epsilon)||u||^2_{W^{1,2}}
    }
    \nonumber
    \\
 &=\frac{1}{\text{Area}(\partial B)} +O(||u||^2_{W^{1,2}}).   
\end{align}
 %\edz{In (6.8) line 1, $O(\epsilon)||u||_{W^{1,2}}$ should be  $O(\epsilon)||u||^2_{W^{1,2}}$.}
Thus
\begin{align}
    h(t)
    &=
    \frac{1}{\text{Area}(M(t))}\bigg( {n \choose k}\text{Area}(M(t))
    +O(||u||_{W^{1,2}}^2)\bigg)
    \nonumber
    \\
    &=
    {n \choose k}
    +O(||u||_{W^{1,2}}^2).
\end{align}
Using the formulas in Lemma 4.1 %\ref{taylor} 
and Lemma 3.1 in \cite{MR2441221},
%\ref{sigmak},
we have $\sigma_k(L)$ equals:
\begin{align}
&\frac{1}{(1+u)^2 + |\nabla u|^2)^{\frac{k+2}{2}}}
 \sum_{m = 0}^{k}
 \frac{(-1)^m{n - m \choose k -m} }{(1+u)^m}
 \bigg((1+u)^2\sigma_m(D^2u)
  \nonumber
\\
&  +
\frac{n+k-2m}{n-m} u^iu_j [T_m]_i^{j}(D^2u)
\bigg )
\nonumber
\\
=&
{n \choose k}(1 -(k+2)u)
(1+2u - \frac{k}{n}\Delta u) + O(\|u\|^2_{W^{1,2}})
%O(u^2, |\nabla u|^2, |D^2u|^2, u \Delta u)
\nonumber
\\
=&
{n \choose k}
(1
- ku -  \frac{k}{n}\Delta u
)
+O(\|u\|^2_{W^{1,2}})
%O(u^2, |\nabla u|^2, |D^2u|^2, u\Delta u).
\end{align}
Thereby
\begin{align}
    G 
    &=
    (-\sigma_k(L) + h)
    =
    k{n \choose k}( u + \frac{1}{n}\Delta u)
     +O(\|u\|^2_{W^{1,2}})
%     O(u^2, |\nabla u|^2, |D^2u|^2, u\Delta u)
    \end{align}
    %\edz{Use $O(\|u\|^2_{W^{1,2}})$ in (6.10), (6.11), (6.12) instead.}
and
\begin{align}
    G ^2
    &=
    k^2 {n \choose k}^2
   \bigg(
   u^2 + \frac{1}{n^2}(\Delta u)^2 + \frac{2}{n}u\Delta u 
   \bigg )
    + O(\epsilon)(\|u\|^2_{W^{1,2}})
%    O(\epsilon)O(u^2, |\nabla u|^2, |D^2u|^2).
\end{align}

Next, as in \eqref{radius_flow_formula}, $u_t = G\sqrt{1 + \frac{|\nabla u|^2}{(1+u)^2}}$. Then, since
$$\sqrt{1 + \frac{|\nabla u|^2}{(1+u)^2}} = 1 + 
O(\|u\|^2_{W^{1,2}})
%O(u^2, |\nabla u|^2),
,$$ we obtain
\begin{align}
    \frac{d}{dt}
    ||u||_{L^2}
    &=2 \intsphere
    u u_t
    \textit{dA}
    \nonumber
    \\
 &=   2 \frac{k}{n} {n \choose k} \intsphere
   ( nu^2 + u\Delta u )\sqrt{1 + \frac{|\nabla u|^2}{(1+u)^2}}
    \textit{dA}
    +
    O(\epsilon)||u||_{W^{2,2}}^2
    \nonumber
       \\
 &=   2 \frac{k}{n} {n \choose k} \intsphere
    nu^2 - |\nabla u|^2
    \textit{dA}
    +
    O(\epsilon)||u||_{W^{2,2}}^2.
\end{align}
Using the formula in Lemma \ref{Bertini_deriv} with $\alpha=1$, we compute
\begin{align}
    \frac{d}{dt}
    \int_{M(t)}
    \sigma_{k-1}\textit{d}\mu_t
    &= 
   - k\intsphere G^2
    (1+u)^n\sqrt{1 + \frac{|\nabla u|^2}{(1+u)^2}}
    \textit{dA}
    \nonumber
    \\
    &=
    - k^3{n \choose k}^2\intsphere
    u^2 + \frac{1}{n^2}(\Delta u)^2
    +\frac{2}{n}u\Delta u
    \textit{dA}
    +
    O(\epsilon)||u||_{W^{2,2}}^2
    \nonumber
     \\
    &=
    - k^3{n \choose k}^2\intsphere
    u^2 + \frac{1}{n^2}(\Delta u)^2
    -\frac{2}{n}|\nabla u|^2
    \textit{dA}
    +
    O(\epsilon)||u||_{W^{2,2}}^2.
\end{align}
%\edz{First line's right hand side of (6.14) should be $(1+u)^n$, not $(1+u)^2$.}

\end{proof}

Now we revisit Proposition \ref{propgen} along the flow, from which we conclude that if $M$ is a nearly spherical surface, where the barycenter is at the origin and $V(\Omega) = V(B)$, then
\begin{align}
\label{prop_ch5_sec2}
    I_{k-1}(\Omega) - I_{k-1}(B)
    &\geq
    {n \choose k-1}
    \frac{k(n-k+1)}{2n}
    ||u||_{L^2}^2
    +O(\epsilon)
     || u||_{W^{1,2}}^2.
\end{align}
This is a weaker statement than Proposition \ref{propgen} because we have left out the gradient term. However, the gradient term is not needed when using the Fraenkel asymmetry in the quantitative quermassintegral inequality.  

In the following theorem, we compare the derivatives of both sides of the inequality \eqref{prop_ch5_sec2}. Importantly, just as in Theorem \ref{main_inversemeancurv_inequality}, we need a condition that ensures the barycenter remains near the origin throughout the flow. We know this condition will be preserved for n-symmetric sets, for example. 
    
For convenience, we recall the statement of Theorem \ref{MeanFlowMainThm}.
\begin{Thm}

Suppose $M(t)$ is a solution of surfaces to the flow \eqref{MeanCurvFlowStatement2}, and at $t_0$ the surface $M(t_0)$ satisfies, for a fixed $K$, that $|\text{bar}(M(t_0))|^2 \leq K\epsilon||u(t_0)||_{W^{2,2}}^2$ and $||u(t_0)||_{W^{2,\infty}} < \epsilon$. Then, for any small $\eta>0$
\begin{align}
\label{meanflow_derivbound}
       \frac{d}{dt}\bigg ( I_{k-1}(\Omega(t_0)) - I_{k-1}(B) \bigg )
       \leq
      (1-\eta) \frac{d}{dt}
     \frac{k(n-k+1)}{2n}
    ||u(t_0)||_{L^2}^2,
\end{align}
where the choice of a sufficiently small $\epsilon>0$ depends on $\eta$ and $K$.
Additionally, if $|\text{bar}(M(t))|^2 \leq K\epsilon||u||_{W^{2,2}}^2$ holds for sufficiently large $t$, then
\begin{align}
\label{liminf_meancurv}
    \liminf_{t\rightarrow \infty}
    \frac{I_{k-1}(\Omega(t)) - I_{k-1}(B)}{
    ||u(t)||_{L^2}^2} &\geq \frac{k(n-k+1)}{2n}.
\end{align}
\end{Thm}
\begin{Rem}
Here we are not able to conclude an analogous statement to Corollary \ref{corollary_inversemean} because we do not have a similar result to Lemma \ref{remain_nearly_spherical}. 
\end{Rem}

\begin{proof}
First, from \eqref{uL2meancurv} $\frac{d}{dt}
     \frac{k(n-k+1)}{2n}
     {n \choose k-1}
    ||u||_{L^2}^2$ equals:
\begin{align}
    &
     \frac{k^2(n-k+1)}{n^2} {n \choose k-1}{n \choose k}
     \intsphere
     nu^2 - |\nabla u|^2
     \textit{dA}
      + O(\epsilon)||u||_{W^{2,2}}^2
     \nonumber
     \\
    =&
     \frac{k^3}{n^2} {n \choose k}^2
     \intsphere
     nu^2 - |\nabla u|^2
     \textit{dA}
      + O(\epsilon)||u||_{W^{2,2}}^2.
\end{align}

From Lemma \ref{poincare_laplace_bound}, when $
|\text{bar}(M(t))|^2 \leq K \epsilon||u||_{W^{2,2}}$, we have  $||\Delta u||_{L^2}^2 \geq 2(n+1) ||\nabla u||_{L^2}^2 +O(\epsilon)||u||_{W^{2,2}}$.  Thus, we can bound the deficit with formula \eqref{sigmakmeancurv},
\begin{align}
    & \frac{d}{dt}(I_{k-1}(M(t)) - I_{k-1}(B))\nonumber\\
     =&
      \frac{k^3}{n^2}{n \choose k}^2\intsphere
  -n^2  u^2 
    -(\Delta u)^2
    +2n
   |\nabla u|^2
    \textit{dA}
    +
    O(\epsilon)||u||_{W^{2,2}}^2
    \nonumber
    \\
 \leq&
  \frac{k^3}{n^2}{n \choose k}^2\intsphere
 - n^2  u^2 
-
   |\nabla u|^2
    \textit{dA}
     + O(\epsilon)||u||_{W^{2,2}}^2
 .
\end{align}

We achieve both \eqref{meanflow_derivbound} and \eqref{liminf_meancurv} by the same reasoning as in Theorem \ref{main_inversemeancurv_inequality}, which is made possible from the convergence of $M(t)$ to a sphere as stated in Lemma \ref{sinestrari_thm}.

\end{proof}

\bibliographystyle{alpha}
\bibliography{references}

\end{document}